\documentclass[final, crop=false, float = true, 11pt]{aims}

\usepackage{paralist}
\usepackage{tiaArticleStyle}

\hypersetup{urlcolor=blue, citecolor=red}

\def\currentvolume{X}
 \def\currentissue{X}
  \def\currentyear{200X}
   \def\currentmonth{XX}
    \def\ppages{X--XX}

\theoremstyle{definition}

\newcommand{\eps}[1]{{#1}_{\varepsilon}}
\renewcommand{\eps} {\varepsilon}
\newcommand{\OID}{OID\xspace}
\usepackage{amsbsy}

\title[Optimal inversion design]
      { Efficient learning methods for large-scale optimal inversion design}

\author[J. Chung, M. Chung, S. Gazzola, and M. Pasha]{}

 \email{jmchung@vt.edu}
 \email{mcchung@vt.edu}
 \email{sg968@bath.ac.uk}
 \email{mpasha3@asu.edu}

\thanks{The first author is supported by NSF grants DMS-1654175 and DMS-1723005.\\
        The second author is supported by NSF grant DMS-1723005.\\
        The third author is supported by EPSRC grant  EP/T001593/1.\\
        The fourth author is supported by NSF grant DMS-1502640.}

\thanks{$^*$ Corresponding author: Matthias Chung, \href{mailto:mcchung@vt.edu}{\tt mcchung@vt.edu}}

\begin{document}
\maketitle

\centerline{\scshape Julianne Chung and Matthias Chung$^*$}
\medskip
{\footnotesize
 \centerline{Department of Mathematics, Academy of Data Science}
   \centerline{Virginia Tech, Blacksburg, VA 24061, USA}
} 

\medskip

\centerline{\scshape Silvia Gazzola}
\medskip
{\footnotesize
 \centerline{Department of Mathematical Sciences}
   \centerline{University of Bath, Bath BA2 7AY, UK}
} 

\medskip

\centerline{\scshape Mirjeta Pasha}
\medskip
{\footnotesize
 \centerline{School of Mathematical and Statistical Sciences}
   \centerline{Arizona State University, Tempe, AZ, USA}
}

\bigskip

 \centerline{(Communicated by the associate editor name)}

\begin{abstract}
In this work, we investigate various approaches that use learning from training data to solve inverse problems, following a bi-level learning approach. We consider a general framework for \emph{optimal inversion design}, where training data can be used to learn optimal regularization parameters, data fidelity terms, and regularizers, thereby resulting in superior variational regularization methods. In particular, we describe methods to learn optimal $p$ and $q$ norms for ${\rm L}^p-{\rm L}^q$ regularization and methods to learn optimal parameters for regularization matrices defined by covariance kernels. We exploit efficient algorithms based on Krylov projection methods for solving the regularized problems, both at training and validation stages, making these methods well-suited for large-scale problems. Our experiments show that the learned regularization methods perform well even when there is some inexactness in the forward operator, resulting in a mixture of model and measurement error.
\end{abstract}

\section{Introduction}\label{sec:Intro}
Inverse problems arise in many important science and engineering applications such as biomedical and astronomical imaging, satellite surveillance, and seismic monitoring \cite{hansen2010discrete,calvetti2007introduction}. Two of the main challenges to solving large-scale inverse problems are (i)~ill-posedness of the problem, whereby small noise or errors in the data can and often do lead to large errors in the solution, and (ii)~the large size of the problem, which for some applications is on the order of millions of observations and billions of unknown parameters. A standard way to solve inverse problems is to follow a variational approach, where solutions are computed by minimizing a pre-determined energy functional that depends upon assumptions regarding the statistical distribution of the observational noise, the forward model, and any prior knowledge about the properties of the unknown solution.
Although a significant amount of research has gone into developing efficient optimization methods to solve variational problems, the formulation of the optimization problem relies on standard assumptions that may not hold in general and that, moreover, may further rely on additional unknown (hyper)parameters.\\

In this work, we describe a general \emph{optimal inversion design} (\OID) framework for solving inverse problems, where the goal is to use available training data to design an optimal energy functional for variational inversion. In order to introduce the \OID\ learning problem, we begin with a discrete linear inverse problem of the form,
\begin{equation}\label{eq:linEq}
    \bfb = \bfA\bfx_{\rm true} +\bfe,
\end{equation}
where $\bfA \in \bbR^{m\times n}$ represents a given forward model that is also known as the parameter-to-observation map, $\bfb\in\bbR^m$ stores available observations corrupted by some unknown additive noise $\bfe\in \bbR^m$, and $\bfx_{\rm true}$
contains unknowns that should be recovered.
We assume that the inverse problem is ill-posed,
and therefore regularization is needed to compute stable, reasonable approximations of $\bfx_{\rm true}$.
The aim of regularization is to incorporate prior knowledge about the solution. There are many forms of regularization ranging from spectral filtering methods to variational regularization methods to iterative regularization, and many combinations and variants of these  \cite{hansen2010discrete,engl1996regularization}.
In its general form, we consider approaches where the regularized solution can be computed as
\begin{equation}\label{eq:inner}
\widehat\bfx(\bftheta) \in \argmin_{\bfx\in\bbR^n}  \ \calJ(\bfx, \bfA, \bfb; \bftheta) + \calR(\bfx; \bftheta),
\end{equation}
where the overall loss is composed of a data fitting term $\calJ$, which incorporates the forward process $\bfA$ and information about the measurement process, such as the noise distribution in the observations $\bfb$, and a regularization functional $\calR$ that integrates prior knowledge of $\bfx_{\rm true}$. While determined by the underlying statistics, the selection of $\calJ$ and $\calR$ is problem dependent and remains a crucial yet heuristic choice for the inversion process \cite{calvetti2007introduction}. Here we assume that such \emph{design} choices may be represented by some \emph{design parameters} $\bftheta \in \bbR^\ell$, often also referred to as \emph{hyperparameters} \cite{dunlop2019hyperparameter}.

Within this work we focus on a particular form of~\eqref{eq:inner} which is given by
\begin{equation}\label{eq:generalLowelevel_}
    \widehat\bfx(\bftheta)  \in \argmin_{\bfx} \ \norm[p]{\bfA\bfx - \bfb}^p + \lambda \norm[q]{\bfL(\bfbeta)\bfx}^q,
\end{equation}
with design parameters $\bftheta = [\lambda; p; q; \bfbeta]$, where $\lambda, p, q \in \bbR^+$ and $\bfbeta \in \bbR^{\ell_\bfbeta}.$
Here, $\norm[s]{\mdot}$ denotes the ${\rm L}^s$-norm for $s\geq 1$ and a homogeneous function without all norm properties for $0<s<1$. This formulation encompasses many popular variational regularization methods. For instance:
\begin{enumerate}
    \item For fixed $\bfL(\bfbeta)=\bfL \in \bbR^{r \times n}$ and $\bftheta = [\lambda,p,q]\t$, problem~\eqref{eq:generalLowelevel_} is an ${\rm L}^p-{\rm L}^q$ type regularized problem,
    \begin{equation}
    \label{eq:lplq}
       \min_\bfx  \ \norm[p]{\bfA \bfx - \bfb}^p + \lambda \norm[q]{\bfL\bfx}^q.
    \end{equation}
 \item For fixed $p =q = 2$ and $\bftheta = [\lambda,\bfbeta]\t$, problem~\eqref{eq:generalLowelevel_} may include a design-dependent operator $\bfL:\bbR^{\ell_{\bfbeta}} \to \bbR^{r \times n}$ in the regularization term, i.e.,
    \begin{equation}\label{eq:parampb}
    \min_\bfx \ \norm[2]{\bfA \bfx - \bfb}^2 + \lambda \norm[2]{\bfL(\bfbeta)\bfx}^2.
    \end{equation}
    Within a Bayesian approach $\bfL(\bfbeta)$ may be regarded as an inverse square root of a positive definite parameteric prior covariance matrix. Consequently, a minimizer of~\eqref{eq:parampb} may then constitute a maximum a posteriori estimate \cite{kaipio1999inverse, kaipio2006statistical}.

\end{enumerate}
Both problems \eqref{eq:lplq} and~\eqref{eq:parampb} depend on the particular choice of the design parameters $\bftheta$, and the main question is how to \textit{optimally} select $\bftheta$?

Assume that we are given the distribution of $\bfx_{\rm true}$ and $\bfe$. Then optimal design parameters $\bftheta$ may be selected by minimizing the Bayes risk, i.e.,
\begin{subequations}
\begin{gather}\label{eq:outer}
		\min_{\bftheta \in \Omega} \ \thf \bbE \norm[2]{\widehat\bfx(\bftheta) -  \bfx_{\rm true}}^2\\
		\text{ subject to } \widehat\bfx(\bftheta) \text{ solving}~\eqref{eq:generalLowelevel_},
\end{gather}
\end{subequations}
where $\Omega$ is a set of feasible design choices and $\bbE$ is the expected value. By minimizing the expected mean squared error \eqref{eq:outer}, the optimal design parameters are expected to perform well on average, leading to reconstructions $\widehat\bfx(\bftheta)$ that minimize the Bayes risk. While other design criteria are available, we focus on this design criterion, which is referred to as \emph{A-design} in the field of optimal experimental design \cite{pukelsheim2006optimal,atkinson1992optimum}.

For problems where the distribution of $\bfx_{\rm true}$ is unknown or not obtainable, but training data are readily available, we consider empirical Bayes risk design problems, where the training data are used to approximate the expected value in \eqref{eq:outer}.  Assume that we are given a set of training data consisting of $J$ true models $\bfx_{\rm true}^1, \ldots, \bfx_{\rm true}^J \in \bbR^n$ and simulated observations $\bfb^1, \dots, \bfb^J \in \bbR^m$, e.g., by data simulation through \eqref{eq:linEq}. Then we consider the empirical Bayes risk \OID\ problem,
\begin{subequations}\label{eq:all_sample}
\begin{gather}
	\widehat\bftheta \in \argmin_{\bftheta \in \Omega}\;\tfrac{1}{2J} \sum_{j = 1}^J \norm[2]{\widehat\bfx^j(\bftheta) -  \bfx_{\rm true}^j}^2 \label{eq:outer_sample}\\
		\text{ subject to}\; \widehat{\bfx}^j(\bftheta) \text{ {solving}}
		~\eqref{eq:generalLowelevel_}
		\text{ for data } \bfb=\bfb^{j}.\label{eq:inner_sample}
\end{gather}
\end{subequations}
In other words, the design problem is a bi-level optimization problem where the goal is to find the parameters $\bftheta$ that minimize the sample average of reconstruction errors for some training set \cite{haber2003learning,haber2008numerical,chung2011designing,calatroni2017bilevel,antil2020bilevel}. The \emph{outer} optimization problem~\eqref{eq:outer_sample} is referred to as the design problem, while the variational regularization problem~\eqref{eq:inner_sample} is referred to as the \emph{inner} problem. \\

\noindent\emph{Overview of main contributions.} In this work, we describe efficient learning techniques to solve the overall design problem~\eqref{eq:all_sample}. Although this framework can incorporate various variational regularization techniques, we focus on the two scenarios described above. Learning the regularization parameter $\lambda$ has been previously considered in various contexts, but to the best of our knowledge, learning optimal values of $p$ and $q$ (i.e., $\bftheta=[\lambda,p,q]$) for the ${\rm L}^p-{\rm L}^q$ regularized problem and learning optimal parameters for covariance kernel matrices (i.e., $\bftheta=[\lambda,\bfbeta]$) have not been considered in an \OID\ framework before.  We will show that these methods can handle various uncertainties in the problem, from mitigating errors in the forward model to resolving unknown parameters in the prior and noise assumptions. Furthermore, we exploit recent developments in Krylov projection methods to efficiently handle the inner problem.

An outline of the paper is as follows. In Section \ref{sec:previouslearning} we provide a brief overview on previous research on learning methods for solving inverse problems. Section \ref{sec:computation} is devoted to computational approaches for learning design parameters in an \OID\ framework including details on iterative projection methods for solving the inner problems \eqref{eq:inner_sample}. In Section \ref{sec:numerics}, we provide numerical results for various image deblurring and tomography applications that demonstrate the effectiveness and benefits of our approaches. Conclusions are provided in Section \ref{sec:conclusions}.

\section{Previous works on learning for inverse problems}
\label{sec:previouslearning}

Supervised learning techniques have gained significant interest in the inverse problems community as a way to combine model-driven and data-driven approaches for solving inverse problems.  A comprehensive overview can be found in \cite{arridge2019solving}.  Two predominant classes of supervised learning approaches have emerged for solving inverse problems: empirical Bayes risk minimization approaches that are related to optimal experimental design techniques \cite{haber2008numerical} and techniques based on deep learning tools such as neural networks and variational autoencoders \cite{lucas2018using}. Supervised training approaches for solving inverse problems were first formally introduced in Haber and Tenorio~\cite{haber2003learning}, where a bi-level optimization problem of the form \eqref{eq:all_sample} was considered for learning optimal parameters for the regularization functional. One of the many advantages of these learning approaches is that the learned (parameterized) regularization functional is tailored to a specific forward operator and noise level of the data. There have been various extensions of this idea (e.g., to learn optimal spectral filters \cite{chung2011designing,chung2017learning}, optimal weighted and multi-parameter Tikhonov parameters \cite{haber2010numerical,holler2018bilevel}, and optimal weighted TV parameters \cite{hintermuller2017optimal}). An additional advantage is that empirical Bayes risk minimization approaches can exploit existing computationally efficient optimization techniques and incorporate a wide range of state constraints \cite{ruthotto2018optimal}. Furthermore, they are general in that different design objective functions can be incorporated \cite{haber2008numerical}, they can be used to learn critical information such as optimal sampling patterns (e.g., for MRI \cite{sherry2020learning}) and design setup for experiments (e.g., for tomography \cite{ruthotto2018optimal}), and they have rich theoretical connections to Bayesian experimental design \cite{alexanderian2016bayesian,huan2014gradient}. The main concerns of this approach include the need to solve an expensive bi-level optimization problem \cite{de2013image,de2017bilevel,calatroni2017bilevel} and bias towards the training set, since reconstructed parameters are only good on average.

The other major class of supervised learning techniques to take root in the inverse problems community consists of methods that exploit deep learning techniques (e.g., neural networks and variational autoencoders), see e.g., review papers \cite{arridge2019solving,ongie2020deep,mccann2017convolutional,lucas2018using}. Initially, deep learning was used mainly for postprocessing of solutions to improve solution quality (e.g., image denoising) or for performing tasks such as classification.  However, deep neural networks are now being considered for solving inverse problems by learning the mapping from observation-to-reconstruction \cite{hammernik2018learning,zhang2017learning} or by learning an appropriate auto-encoder network (e.g., a generative adversarial network) to serve as a proxy for the regularizer \cite{lunz2018adversarial,haltmeier2020regularization,li2020nett,prost2021learning}. However, a major disadvantage for many of these network learning approaches is that due to the large number of network parameters, a massive amount of training data is required, which may not be readily available. Furthermore, it is important to have a well-tuned network and a good choice of parameters (e.g., batch size, epochs, learning rate) for the stochastic optimization algorithms, prior to training the networks.

Although there have been significant developments in both classes of supervised learning approaches for inverse problems, there are still open problems. In particular, as mentioned in Section \ref{sec:Intro}, we are interested in learning the appropriate ${\rm L}^p$ and ${\rm L}^q$-norms along with the regularization parameter $\lambda$ in \eqref{eq:lplq}. This problem is most related to the work by De los Reyes and Schönlieb \cite{de2013image}, where parameter learning methods were considered for learning the noise model in variational image denoising by estimating weights for different noise models. However, rather than consider a weighted version of pre-determined noise models, our approach seeks an appropriate ${\rm L}^p$-norm to resolve any errors or uncertainties in the data-fit term and combines it with an optimally-selected ${\rm L}^q$-norm for the regularization term. For the problem of learning optimal parameters for a regularization operator \eqref{eq:parampb}, special cases have been considered in supervised learning frameworks (e.g., \cite{haber2003learning} considered different regularization terms for different regions of the solution and \cite{antil2020bilevel} considered a bilevel optimization learning framework for learning the fractional Laplacian parameter).  However, learning the kernel parameters in an \OID\ framework remains a challenge, especially when the prior precision matrix (i.e., the inverse of the covariance matrix) or its square root $\bfL(\bfbeta)$ is not readily available for all design parameters $\bfbeta$. In general, significant computational challenges may arise within large-scale bilevel optimization problem~\eqref{eq:all_sample}, which we address next.

\section{Computational \OID\ for variational inverse problems}
\label{sec:computation}

In this section, we describe computational approaches for the \OID\ problem,
\begin{subequations}\label{eq:oid}
\begin{gather}
	\widehat\bftheta \in \argmin_{\bftheta \in \Omega} \ \calP(\bftheta) =  \tfrac{1}{2J} \sum_{j = 1}^J \norm[2]{\widehat\bfx^j(\bftheta) -  \bfx_{\rm true}^j}^2 \label{eq:oid_outer}\\
		\text{ s.t. }\; \widehat\bfx^j(\bftheta) = \argmin_{\bfx} \ \norm[p]{\bfA\bfx - \bfb^j}^p + \lambda \norm[q]{\bfL(\bfbeta)\bfx}^q \quad j = 1,\ldots, J,\label{eq:oid_inner}
\end{gather}
\end{subequations}
where $\bftheta = [\lambda, p, q, \bfbeta]$ with $\lambda, p, q >0$.

Bi-level optimization problems such as \eqref{eq:oid} are notoriously difficult to solve. For instance, simple non-convexity in the inner problem (such as those encountered when $p,q<1$ in \eqref{eq:oid_inner}) may lead to discontinuities in the outer design problem,  \cite{dempe2002foundations,sinha2017review}. For example, consider the toy problem $\min_\theta \ \widehat x(\theta)$ where the inner problem $\widehat x(\theta) = \argmin_x \ (x-1)^2 (x+1)^2 + \theta x$ is non-convex in $x$. For $\theta = 0$ two global minima $\widehat x(0) = \pm 1$ exist and for any $\theta<0$ and $\theta>0$ we have $\widehat x(\theta) >1$ and $\widehat x(\theta) <-1$, respectively. Hence the outer design function is discontinuous at $\theta = 0$.

Various approaches exist to address bi-level optimization problems. One approach is to cast the inner problem as a constraint and utilize ``off-the-shelf'' constrained optimization methods, such as augmented Lagrangians or interior-point methods. Computational challenges arise in this approach, since the inner problem results in non-standard equality and inequality constraints \cite{calatroni2017bilevel,de2017bilevel,antil2020bilevel,holler2018bilevel}. Another approach commonly used in the PDE constrained optimization literature is to eliminate the constraints by approximately solving for $\widehat{\bfx}^j$, yielding a reduced problem. Such approaches were used to compute optimal error filters in \cite{chung2011designing,haber2008numerical,haber2010numerical}; however, for regularized solutions that have a nontrivial dependence on $\bftheta$ (e.g., for general variational regularization methods~\eqref{eq:inner} where the inner problem does not admit a closed form solution) such methods do not apply.  In a third approach, the potentially discontinuous outer design problem is treated with non-gradient based global optimization methods such as evolutionary methods \cite{weise2009global,sinha2017review}. Note that global optimization methods tend to be significantly more expensive than methods from convex optimization. Hence, special care must be taken in reducing the overall computational cost.

We approach the computational bottleneck from two directions. First, since the inner problem consists of a variational linear inverse problem, we take advantage of recently developed and highly efficient iterative solvers. Specifically, for the ${\rm L}^p-{\rm L}^q$ type problems, we use a majorization-minimization (MM) approach together with generalized Krylov subspaces (GKS), dubbed MM-GKS \cite{buccini2020modulus}; for the parametric kernel learning problem, we use the generalized Golub-Kahan (genGK)-based method, sometimes in its hybrid version, dubbed genHyBR \cite{chung2017generalized}. We refer to the discussions in Sections~\ref{sub:learningexponents} and \ref{sub:operators}, respectively, for more details. Second, we utilize surrogate optimization techniques, also referred to as Bayesian optimization, for the outer problem~\eqref{eq:oid_outer}; see \cite{gramacy2020surrogates,osborne2009gaussian}. The advantage of surrogate optimization methods is that they construct a surrogate objective function and evaluate the surrogate instead of the true objective to find global minimizers\footnote{In Bayesian optimization, it is common to consider equivalent maximization problems.}, thereby reducing the overall number of inner solves~\eqref{eq:oid_inner}.

More precisely, a surrogate optimization method takes samples of the objective function given as $\calS_K = \{\bftheta_k, \calP(\bftheta_k)\}_{k = 1}^K$ and builds a surrogate model $s_{K}:\Omega \to \bbR$ by extrapolating the objective function~\eqref{eq:oid_outer} beyond the sample points $\calS_K$. For instance, surrogate models may be constructed using radial basis functions \cite{urquhart2020surrogate} or Gaussian processes \cite{gramacy2020surrogates}. Typically the surrogate model matches $\calP$ exactly at points $\bftheta_k$, $k = 1,\ldots, K$, hence interpolating the true objective function $\calP$ at $\bftheta_k$, $k = 1,\ldots, K$. From the surrogate model $s_K$, a merit or so called \emph{acquisition function} $m_K:\Omega \to \bbR$ is constructed that balances the trade-off between exploitation and exploration~\cite{archetti2019acquisition}. A commonly used acquisition function is the expected improvement function, where the surrogate model predicts low objective function values by means of known sample locations and values, as well as taking into account uncertainty of unexplored regions. In this work we utilize standard Matlab libraries for surrogate optimization provided by the global optimization toolbox \cite{Matlabsurrogateopt} for outer problem \eqref{eq:oid_outer}. Next we describe two computational OID problems: learning optimal $p$ and $q$ norms and learning optimal hyperparameters for the prior.

\subsection{Learning optimal $p$ and $q$ norms} {\bf \OID\ with $\bftheta = [\lambda, p, q]$.}

\label{sub:learningexponents}
Consider OID problem \eqref{eq:oid} where $\bfL(\bfbeta) = \bfL$ is fixed and $\bftheta = [\lambda, p, q]$, i.e., learning the optimal regularization parameter, data fidelity norm, and regularization norm.

There are various reasons why one would want to learn an optimal $p.$ While it is well-known that $p=2$ should be considered when the noise follows an i.i.d.~Gaussian distribution and that $0<p<2$ should be considered when the available data are perturbed by non-Gaussian noise, it is unclear what to use when there is a mixture of noise corrupting the data.  For specific statistical models of noise, e.g., mixed Gaussian and Poisson noise that arise from Charge Coupled Device detectors, a reformulation to a weighted least-squares problem has been considered, see e.g., \cite{bardsley2006covariance,calatroni2017infimal, kubinova2018robust, buccini2020modulus} and references therein. However, the reformulation relies on an approximation using knowledge about the noise statistics, which is not necessarily available in practice.

Moreover, learning $p$ can be relevant when the forward operator $\bfA$ used to solve the inner problem (\ref{eq:oid_inner}) is inexact.  That is, the adopted forward model is (slightly) different from the one used to generate the training data; e.g., deblurring problems using erroneous point spread functions or tomographic reconstruction problems with slightly mismatched projection angles.  Estimating and correcting for model errors represent important yet challenging tasks when solving inverse problems. For problems where the user has strong knowledge about the parameterization of the forward model, there are sophisticated ways of accounting for inexactness in the forward operator, see e.g., \cite{chung2010efficient,riis2021computed}. In a learning context, recent approaches to learn implicit and explicit corrections to the operator using neural networks was considered in \cite{lunz2021learned} and learning non-Gaussian models was considered in \cite{smyl2021learning}. However, for many scenarios where a good parameterization does not exist or the goal is not necessarily to determine the model correction itself, we show that it is possible to mitigate inexactness in the forward model as well as resolve any faults in the noise assumptions by determining a better norm for the data-fit term.  That is, we use the OID framework to determine a proper choice of the norm in the data-fidelity function, which is purely informed by the availability of training data, to mitigate any effects of inexactness in the forward operator or faults in our noise assumptions.

Learning an optimal value for $q$ in the regularization term is important as well, as this encodes prior knowledge about the solution.  The most common choice is Tikhonov regularization ($q=2$), but for promoting sparsity in the solution, $q=1$ provides a numerically appealing approximation to the computationally NP-hard $q=0$ case.
More recently, regularization techniques that allow a generic choice of $q>0$ have been developed \cite{huang2017majorization, chung2019}. Nevertheless, for such techniques, the choice of a suitable $q$ that accommodates the properties of the desired solution is not always obvious, hence learning $q$ becomes crucial in many applications where its choice can be informed by training data.

Thus, we consider OID problem \eqref{eq:oid} with $\bftheta = [\lambda, p, q]$, where the efficiency of the approach relies on the ability to quickly and accurately compute ${\rm L}^p-{\rm L}^q$ regularized solutions (i.e., solving \eqref{eq:oid_inner} with $\bfL(\bfbeta)=\bfL$ fixed).  This can be challenging, especially in a large-scale setting.  Although various optimization methods such as primal-dual gradient descent methods  \cite{zhu2008efficient, esser2010general, chambolle2011first} could be used, we consider iterative projection methods, which approximate $\widehat{\bfx}(\bftheta)$ by solving \eqref{eq:oid_inner} in a reduced-dimensional (projected) subspace; such methods can be considered as special instances of MM strategies \cite{hunter2004tutorial,lange2016mm}, as summarized below.

Rewriting the $s$-norm as $\norm[s]{\bfx} = \left(\sum_{j = 1}^n |x_j|^{s-2}\, x_j^2\right)^{1/s}$ and by approximating $|x| \approx (x^2+\eps^2)^{1/2} =: \phi_{\eps}(x)$ with $\eps>0$  for $0<s\leq 1$ and $\epsilon =0$ for $s>1$, we obtain (an approximation of) $\norm[s]{\bfx}^s$ by
\begin{equation}
\norm[s]{\bfx}^s \approx  \sum_{j = 1}^n
(x_j^2+\eps^2)^{(s-2)/2}\  x_j^2
= \sum_{j = 1}^n \phi_\eps(x_j)^{s-2}\,  x_j^2.
\end{equation}
By defining $\bfS_{s,\eps}(\bfx)$ as a diagonal matrix dependent on $\bfx$, with
\begin{equation}
    \bfS_{s,\eps}(\bfx) = \diag{\left[\phi_\eps(x_1)^{s-2},\ldots,\phi_\eps(x_n)^{s-2}\right]}\,,
\end{equation}
we get
\begin{equation}
\norm[s]{\bfx}^s \approx \norm[2]{(\bfS_{s,\eps}(\bfx))^{1/2}\ \bfx}^2,
\end{equation}
where we define the square root elementwise. Hence,
\begin{equation}\label{eq:smoothlplq}
     \norm[2]{(\bfS_{p,\eps}(\bfA\bfx-\bfb))^{1/2}\ (\bfA\bfx-\bfb)}^2 + \lambda \norm[2]{(\bfS_{q,\eps}(\bfL\bfx))^{1/2}\ \bfL\bfx}^2
\end{equation}
is a sufficiently smooth approximation of the objective function in~\eqref{eq:lplq}. Assuming an approximation $\bfx_k$ to $\bfx_{\rm true}$ is available, we consider the  \emph{quadratic tangent majorant} of~\eqref{eq:smoothlplq} at $\bfx_k$ (omitting a constant term), i.e.,
\begin{equation}\label{eq: QuadraticMajorantQ}
\calM \left(\bfx, \bfx_k\right) = \norm[2]{\left(\bfS_{p,\eps}^k\right)^{1/2}(\bfA\bfx-\bfb)}^{2}
+ \lambda \norm[2]{\left(\bfS_{q,\eps}^k\right)^{1/2}\bfL\bfx}^{2},
\end{equation}
where we defined
\begin{equation}\label{genWeights}
    \bfS_{p,\eps}^k = \bfS_{p,\eps}(\bfA\bfx_k-\bfb) \qquad \mbox{and} \qquad \bfS_{q,\eps}^k = \bfS_{q,\eps}(\bfL\bfx_k)\,;
\end{equation}
for details see \cite{huang2017majorization}. Given a point $\bfx_k$, we compute $\bfx_{k+1}$ as an approximate solution minimizing~\eqref{eq: QuadraticMajorantQ}. This process is iterated to approximate a solution of~\eqref{eq:lplq}, and it is referred to as MM.

Classical methods for MM (which, in this particular instance, coincide with IRLS methods \cite{LS}) involve minimizing \eqref{eq: QuadraticMajorantQ} by, e.g., applying CGLS, and result in time-consuming inner-outer iterative strategies. Recently developed strategies bypass classical IRLS schemes and approximate a solution of minimizing \eqref{eq:smoothlplq} by simultaneously computing a new approximation $\bfx_{k+1}$ and updating the weights $\bfS_{p,\eps}^{k+1}, \bfS_{q,\eps}^{k+1}$ in \eqref{genWeights}. These methods involve projections on generalized Krylov subspaces (GKS) \cite{huang2017majorization, buccini2020modulus}, and we refer to them as MM-GKS.

Specifically, the GKS-based solver considered here computes $\widehat{\bfx}(\bftheta)$ starting from an initial approximate solution $\bfx_{0}$ belonging to an initial approximation subspace ${\rm ran}(\bfV_{0}^{\rm GKS})={\rm ran}(\bfV_{h})$ generated by, e.g., performing $1\leq h\ll\min\{m,n\}$ steps of Golub--Kahan bidiagonalization applied to $\bfA$ with initial vector $\bfb$. Then, at the $(k+1)$st iteration, one computes the (skinny) QR factorizations,
\begin{equation}\label{eq: skinnyQR}
\bfS_{p,\eps}^{k} \bfA \bfV_{k+1}^{\rm GKS} = \bfQ_{p} \bfR_{p}, \quad
\bfS_{q,\eps}^{k} \bfL \bfV_{k+1}^{\rm GKS} = \bfQ_{q} \bfR_{q}.
\end{equation}
where $\bfV_{k+1}^{\rm GKS}=[\bfV_{k}^{\rm GKS}, \bfv_{\rm new}]$ and $\bfv_{\rm new}$ is the normalized residual vector $\bfA\t (\bfA{\bfx}_k - \bfb) + \lambda\bfL\t\bfL{\bfx}_k$.
The $(k+1)$st approximate solution reads
${\bfx}_{k+1}=\bfV_{k+1}^{\rm GKS}{\bfy}_{k+1}\in{\rm ran}(\bfV_{k+1}^{\rm GKS}),$ where
\begin{equation}\label{eq: normalEq_QR}
{\bfy}_{k+1}=\argmin_{\bfy\in\bbR^{k+1}}\ \norm[2]{\bfR_{p}\bfy-\bfQ\t_{p}\left(\bfS_{p,\eps}^{k}\right)^{1/2}\bfb}^2 + \lambda\|\bfR_{q}\bfy\|_2^2,
\end{equation}
and where the projected problem is obtained by
 plugging in the factorizations in \eqref{eq: skinnyQR} into the functional \eqref{eq: QuadraticMajorantQ}.
GKS-based solvers can be applied to many instances of \eqref{eq:lplq}, provided that matrix-vector products with $\bfL$ are cheap to compute, and $k\ll\min\{m,n\}$.

\subsection{Learning design-dependent operators} {\bf \OID\ with $\bftheta = [\lambda; \bfbeta]$.}
\label{sub:operators}
Learning approaches can also be used to estimate hyperparameters for regularization functionals that belong to a parametric family of regularizers (e.g., those defined from a kernel function). We consider OID problem \eqref{eq:oid} where $p=q=2$ and $\bfL(\bfbeta)$ and its inverse are not readily available, but matrix vector multiplications with $\bfQ(\bfbeta) = (\bfL(\bfbeta)\t \bfL(\bfbeta))^{-1}$ can be done efficiently. For example, with Gaussian random fields, the entries of the prior covariance matrix are computed directly as $\bfQ_{ij}(\bfbeta) = \kappa(r_{ij}; \bfbeta)$ where $\kappa(\,\cdot\,;\bfbeta)$ is a covariance kernel function that depends on some parameters in $\bfbeta$ and $r_{ij} = \norm[2]{\bfz_i- \bfz_j}$, with $\bfz_i$ corresponding to spatial points in the domain. Although the matrix $\bfQ(\bfbeta)$ may be dense and the inverse or symmetric factorization is not available, matrix-vector multiplications with $\bfQ(\bfbeta)$ can often be done efficiently.

We consider two families of covariance matrices that are built from parameterized kernels: the squared exponential covariance matrix and the Mat$\acute{\text{e}}$rn covariance matrix \cite{williams2006gaussian}.
Given a hyperparameter $\beta$ that plays the role of the characteristic length-scale, the \emph{squared exponential kernel} is defined as
\begin{equation}
\label{eq:sqexp}
\kappa(r; \beta) = \text{exp}\left(-\frac{r^2}{2\beta^2}\right).
\end{equation}
Given two hyperparameters $\beta_1$ and $\beta_2$ that define the smoothness and length scale respectively, the \emph{Mat$\acute{\text{e}}$rn kernel} is defined as
\begin{equation}
\label{eq:matern}
\kappa(r; \beta_1, \beta_2) =  \frac{1}{2^{\beta_1-1}\Gamma(\beta_1)} \left( \frac{\sqrt{2 \beta_1}r}{\beta_2}\right)^{\beta_1}K_{\beta_1}\left(\frac{\sqrt{2 \beta_1} r}{\beta_2}\right),
\end{equation}
where $\Gamma(\,\cdot\,)$ is the Gamma function and $K_{\beta_1}(\,\cdot\,)$ is the modified Bessel function of the second kind of order $\beta_1$. Note that, oftentimes in the literature, the Mat$\acute{\text{e}}$rn parameters are denoted as $\nu = \beta_1$ and $\ell = \beta_2,$ where simplifications of the kernel function can be made for half integers $\nu =p+1/2, p \in \bbN^+.$  We do not impose this constraint here.

In most inverse problems settings, the kernel parameters must be selected \textit{prior} to solving the inverse problem and oftentimes appropriate choices come from expert knowledge.  There exist various approaches in Bayesian statistics for estimating hyperparameters for covariance functions (e.g., cross-validation and maximum likelihood) \cite{williams2006gaussian}.  The process, which is referred to as model selection, seeks to estimate the hyperparameters directly from the data, but these methods can be computationally infeasible, especially for large-scale problems.  For inverse problems in imaging, semivariogram methods were considered in \cite{brown2020semivariogram} for estimating Mat$\acute{\text{e}}$rn parameters, but this approach only works for problems where the observation grid and the solution grid are the same (e.g., in deblurring and denoising). We remark that learning approaches that use training data to estimate parameters defining the regularizer have been considered in \cite{haber2003learning,antil2020bilevel,chung2011designing}; however, contrary to existing methods that work with the precision matrix directly, here we consider regularizers that arise in Bayesian approaches and that correspond to prior covariance matrices defined using parametric kernel functions.

We exploit genGK approaches for efficient inner solves \eqref{eq:oid_inner} requiring only matrix-vector products with the prior covariance matrix. More specifically, we are interested in solving \eqref{eq:linEq} where $\bfe \sim \mathcal{N}(\bfzero, \bfI)$ and $\bfx \sim \mathcal{N}(\bfzero, \lambda^{-1}\bfQ(\bfbeta))$ where $\bfQ(\bfbeta)$ is defined above. By Bayes' formula,
\begin{equation*}
\pi(\bfx|\bfb)\propto \pi(\bfb|\bfx)\pi(\bfx) \propto \exp\left(-\|\bfA\bfx - \bfb\|_2^2 - \lambda\bfx\t\bfQ(\bfbeta)^{-1}\bfx\right).
\end{equation*}
The maximum a posteriori approximation of $\bfx$ can be found by minimizing the negative log-likelihood of $\pi(\bfx|\bfb)$, i.e.,
\begin{equation}\label{eq: log_likelihood}
\widehat\bfx(\bftheta)
= \argmin_{\bfx}{\|\bfA\bfx - \bfb \|_2^2 + \lambda\|\bfx\|_{\bfQ(\bfbeta)^{-1}}^2}
\end{equation}
with $\bftheta = [\lambda;\bfbeta]$ which, since $\bfQ(\bfbeta) = (\bfL(\bfbeta)\t \bfL(\bfbeta))^{-1}$, is equivalent to~\eqref{eq:parampb}.
In this setting, an iterative projection method based on the genGK bidiagonalization can be used to approximate \eqref{eq: log_likelihood}. After performing a change of variables (to avoid computations with $\bfQ(\bfbeta)^{-1}$), the $k$-th iteration of the genGK method is given by
\begin{gather*}
\widehat{\bfx}_k(\bftheta)=\bfQ(\bfbeta)\bfV_k^{\rm genGK}\widehat{\bfy}_k(\bftheta),\\
\mbox{where }\;\widehat{\bfy}_k(\bftheta)=\argmin_{\bfy\in\bbR^k} \ \norm[2]{\bfB_k^{\rm genGK}\bfy-\|\bfb\| \bfe_1}^2 + \lambda\|\bfy\|_2^2\,.
\end{gather*}
The matrices above satisfy the partial genGK matrix factorization, i.e.,
\begin{gather*}
\bfA\bfQ(\bfbeta)\bfV_k^{\rm genGK}=\bfU_{k+1}^{\rm genGK}\bfB_k^{\rm genGK},\\
\mbox{with $\bfV_k^{\rm genGK}\in\bbR^{n\times k}$, \quad \mbox{and} \quad$\bfB_k^{\rm genGK}\in\bbR^{(k+1)\times k}$},
\end{gather*}
together with another similar factorization involving $\bfA\t$.
We refer to \cite{chung2017generalized} for the original derivation.

We conclude this section by mentioning that an upside of all the solvers for \eqref{eq:oid_inner} described so far is that $\lambda$ can be adaptively set during the iterations, i.e., they can be reformulated as so-called ``hybrid methods''. When solving \eqref{eq:oid} where $\lambda$ is a design parameter that is fixed for each instance of the inner problem \eqref{eq:generalLowelevel_}, we will not take advantage of this feature of hybrid methods.  However, we may still be able to exploit this feature of hybrid  methods, for OID where $\bftheta=\bfbeta$, and $\lambda$ is selected automatically.  Numerical comparisons will be presented in Section~\ref{sub:numerical_kernel}.

\section{Numerical experiments}
\label{sec:numerics}
In this section, we provide \OID\ examples to show that learned regularization methods perform well for various inverse problems. In Section~\ref{sub:numerical_pq}, we consider an example from image deblurring, where we learn optimal norms for both the data fit and the regularization term, in addition to an optimal regularization parameter, in order to handle different noise types and to mitigate impacts from an imprecise forward operator.  Then, in Section~\ref{sub:numerical_kernel}, we consider an example from tomographic reconstruction, where optimal parameters are found for parametric prior covariance matrices. For all of the experiments, we assess the quality of a reconstructed solution using the Relative Reconstruction Error (RRE) norm defined by ${\rm RRE}(\bfx)=\frac{\|\bfx-\bfx_{\rm true}\|_2}{\|\bfx_{\rm true}\|_2},$ for some reconstruction $\bfx.$

\subsection{\OID\ with $\pmb{\theta}=[\lambda,p,q]\t$}
\label{sub:numerical_pq}
\begin{figure}
    \centering
    \begin{tabular}{cccc}
    \multicolumn{4}{c}{training images}\\
    \includegraphics[width=.2\textwidth]{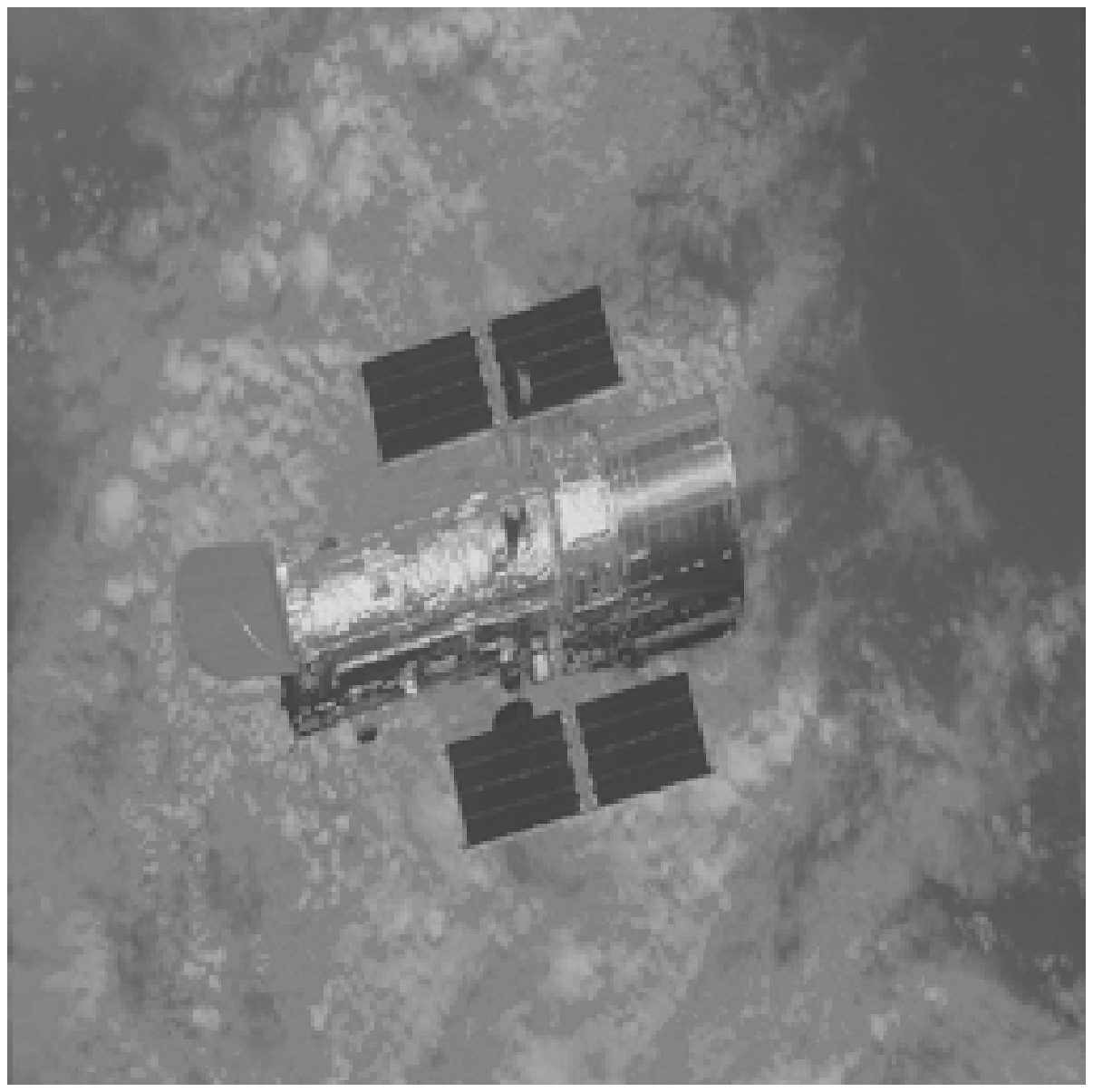} &
    \includegraphics[width=.2\textwidth]{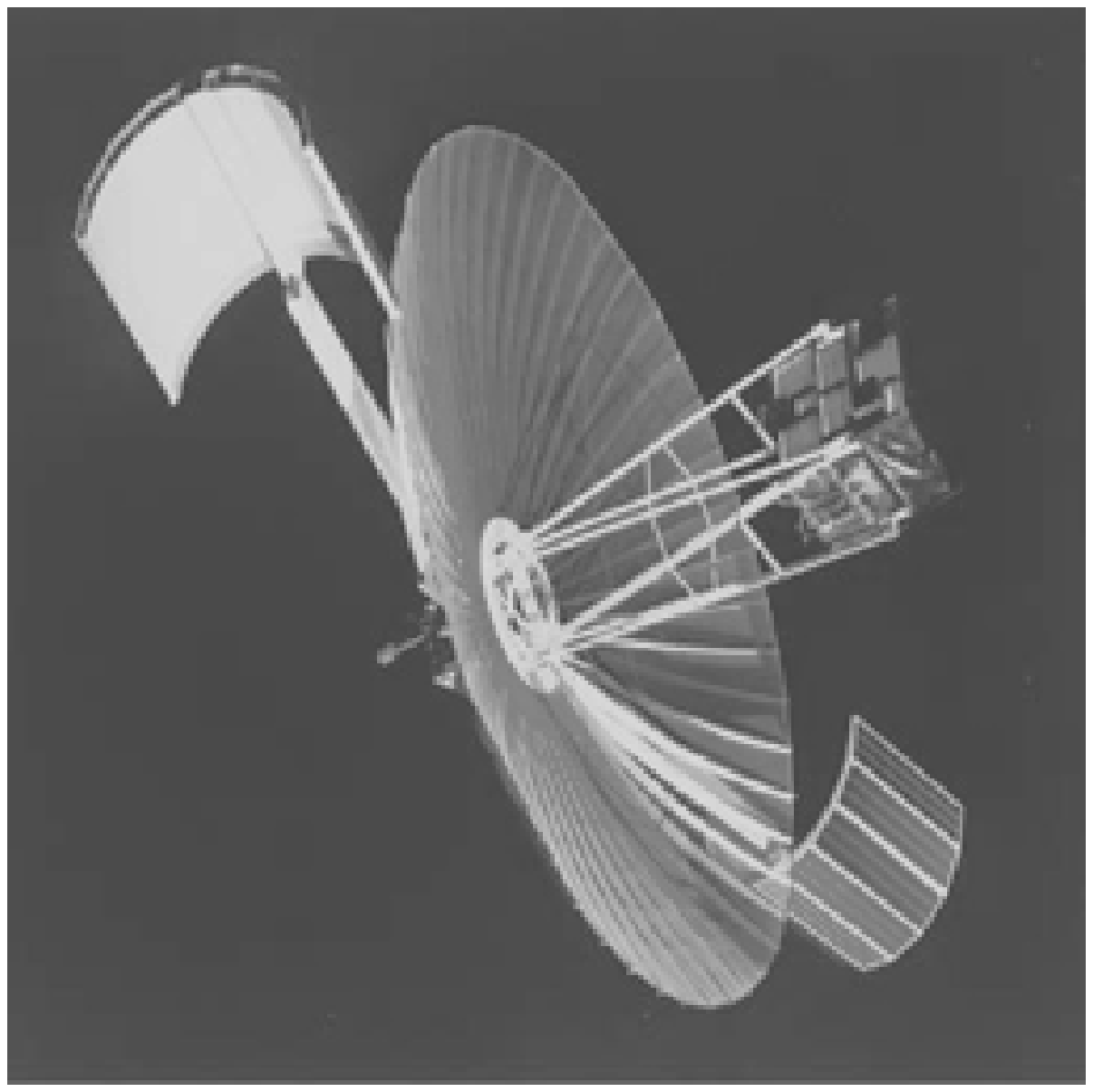} &
    \includegraphics[width=.2\textwidth]{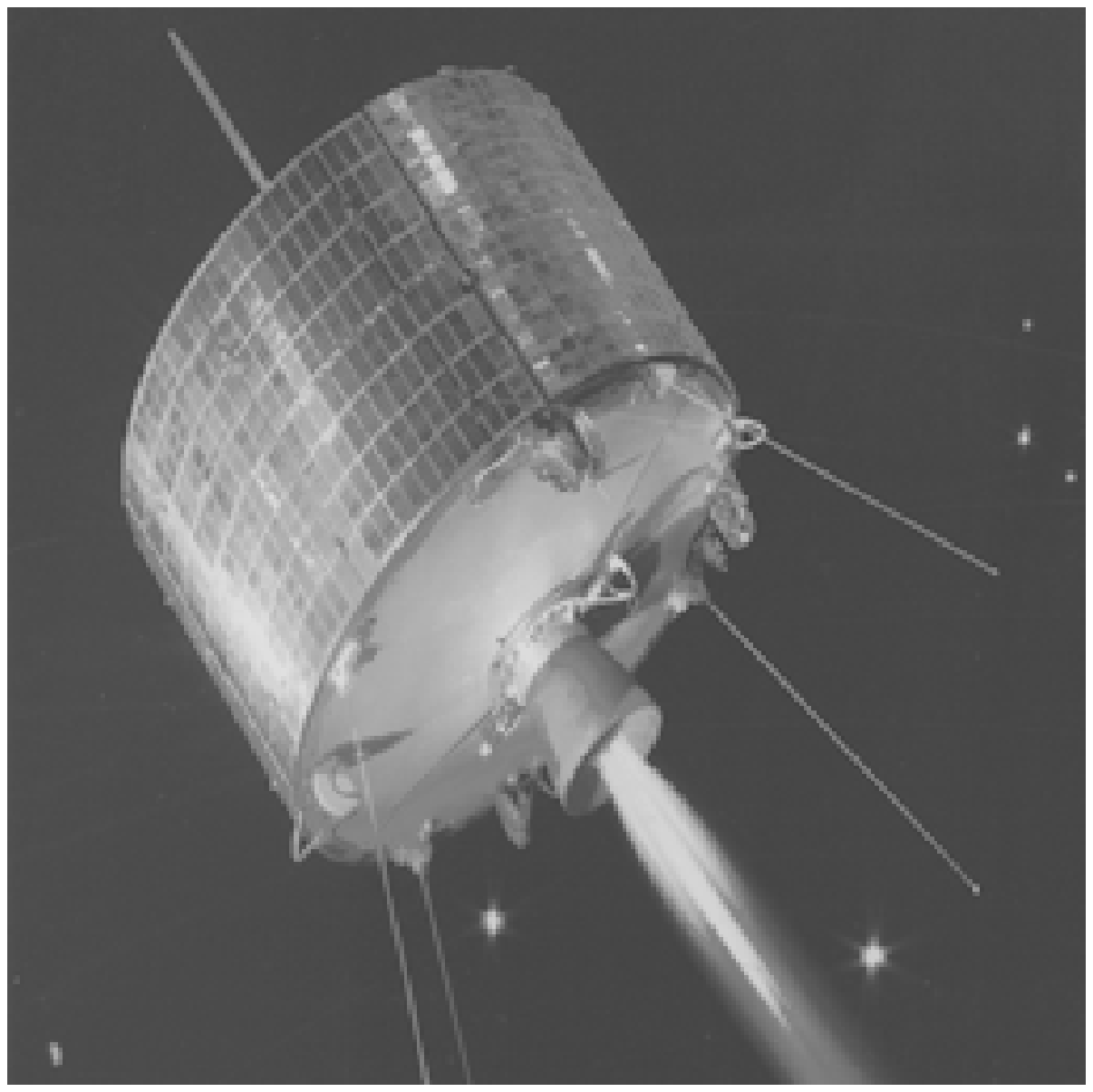} &
    \includegraphics[width=.2\textwidth]{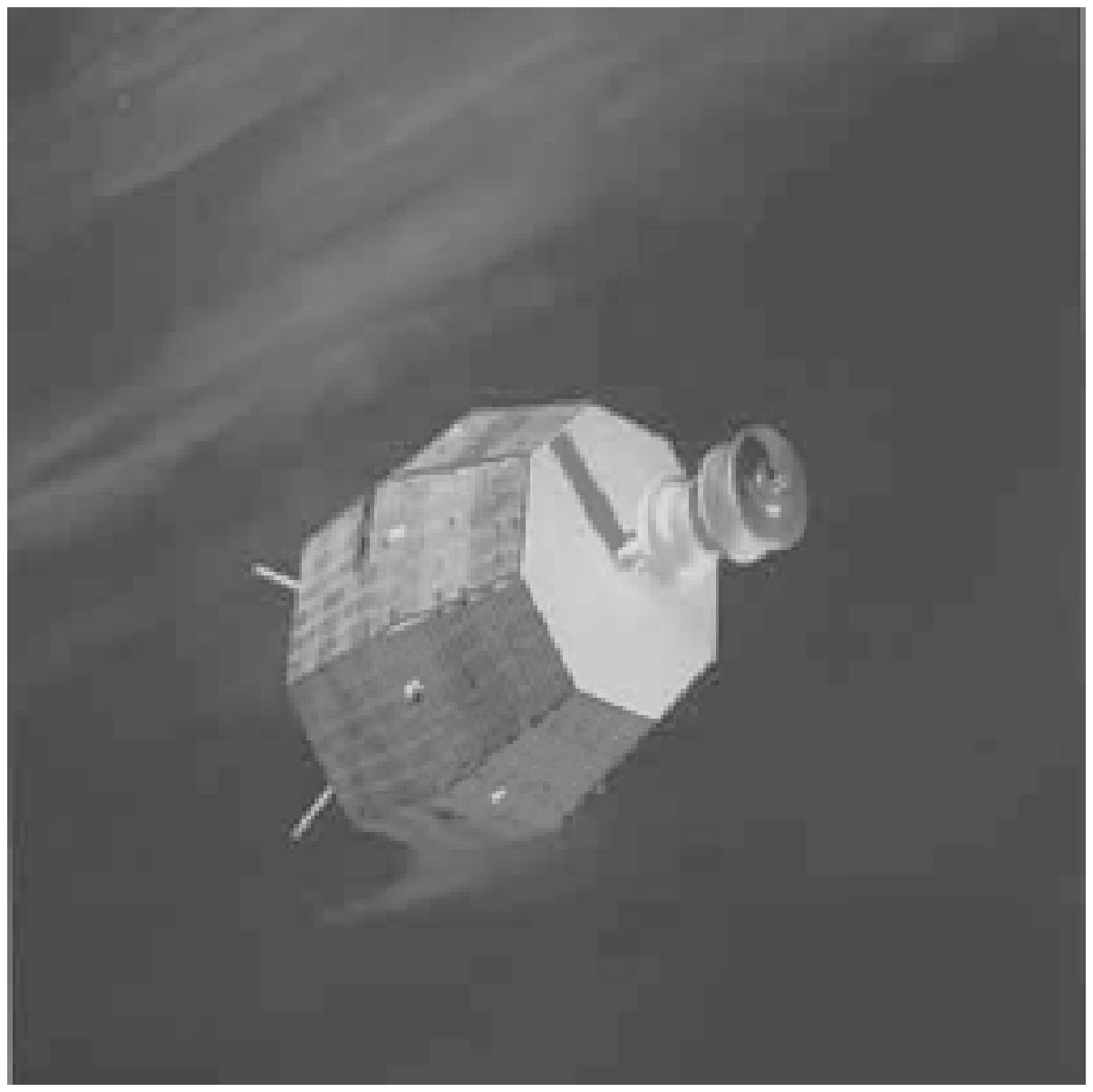}\\
    \multicolumn{4}{c}{validation images}\\
    \includegraphics[width=.2\textwidth]{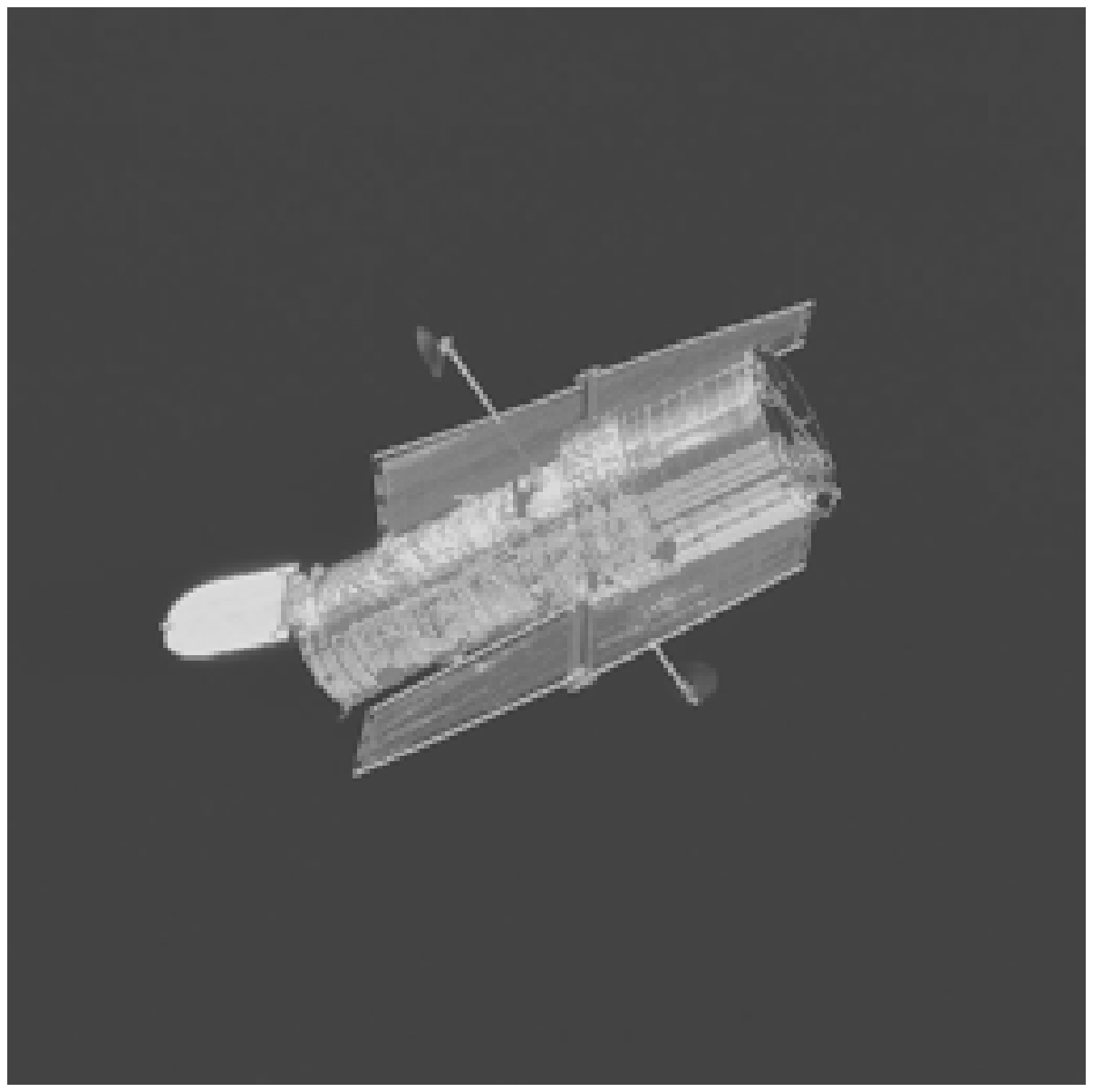} &
    \includegraphics[width=.2\textwidth]{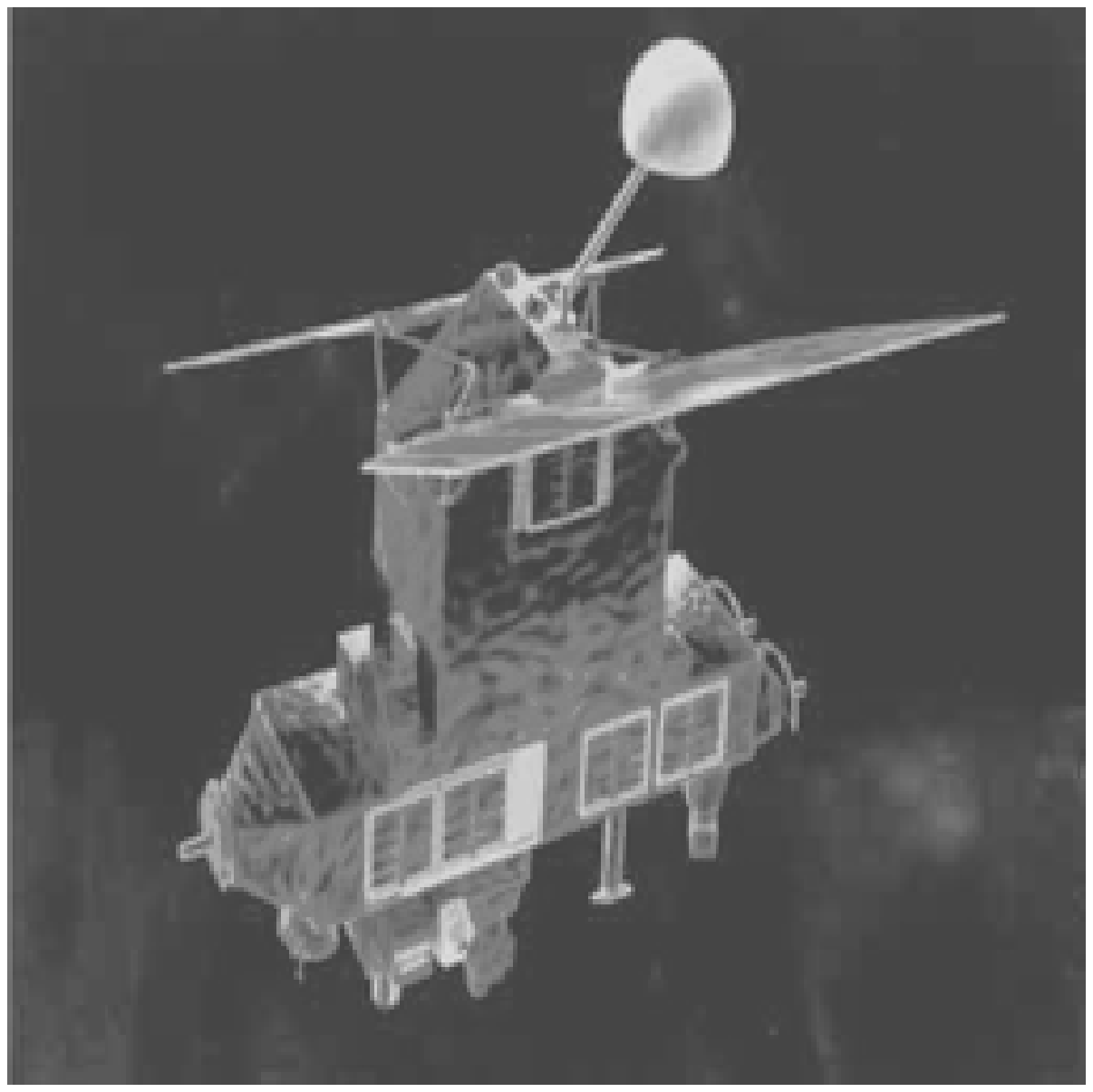} &
    \includegraphics[width=.2\textwidth]{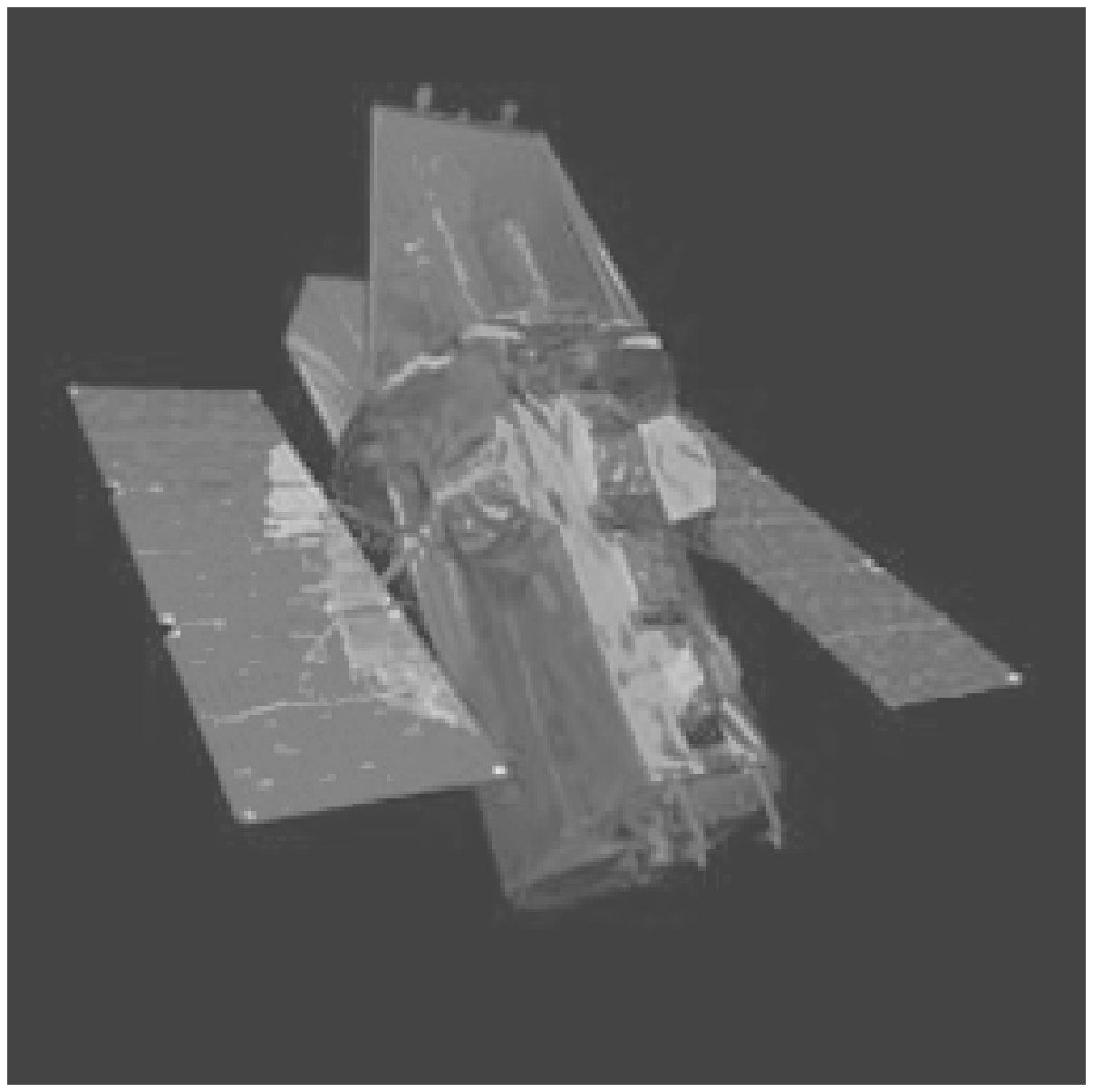} &
    \includegraphics[width=.2\textwidth]{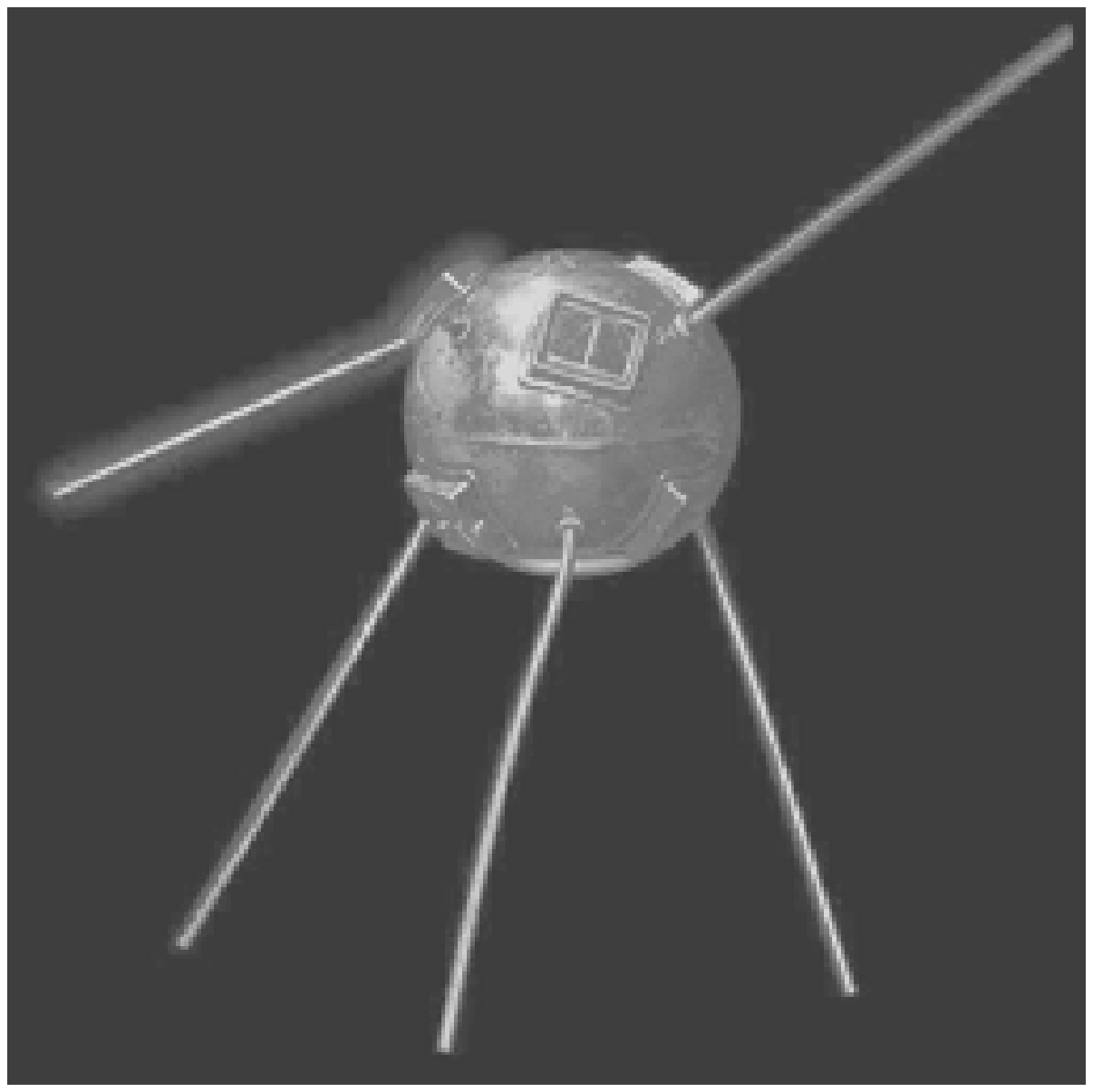}
    \end{tabular}
    \caption{Four prototype true images used for generating the training set (top row) and validation set (bottom row) in the OID experiment with $\bftheta=[\lambda,p,q]\t$.}
    \label{fig:trainvalidsat}
\end{figure}

The goal of this section is to investigate the performance of \OID\ for learning optimal parameters $\lambda, p$, $q$, with $\bfL=\bfI$,  for image deblurring.  For the training and validation datasets, we consider satellite images obtained from the NASA website \cite{NASA}, where each image contains $256\times 256$ pixels. We use 10 images of satellites with 8 random affine transformations, giving a total of 80 training images, and 5 images of satellites with 6 random affine transformations, giving 30 validation images. Samples of the training and validation images are provided in Figure~\ref{fig:trainvalidsat}.

For the forward model, we consider a blurring process defined by an isotropic Gaussian blur centered at location $(\chi_1,\chi_2)$, where the point spread function $\bfP$ has entries
\begin{equation}\label{GaussPSF}
[P]_{i,j}={c(\sigma_1,\sigma_2)}\exp
\left(-\frac{(i - \chi_1)^2}{2\sigma_1^2}
-\frac{(j - \chi_2)^2}{2\sigma_2^2}
\right)\,,
\end{equation}
where $c(\sigma_1,\sigma_2)$ is a scaling factor. In the following we use the notation $\bfA=\bfA(\sigma_1,\sigma_2)$ to highlight the dependence of the matrix $\bfA$ on the blurring parameters, and we consider periodic boundary conditions.

The observed image is obtained as in \eqref{eq:linEq}, with $\bfA(2.5,2.5)$ and $\bfe$ being impulse noise, with noise level selected uniformly at random between 10\% and 50\%. More specifically impulse noise is obtained when the entries of the vector $\bfb$ are constructed as follows
\begin{equation*}
\bfb_i=\left\{\begin{array}{ll}
(\bfA \bfx_{\rm true})_i&\mbox{with probability $1-\eta$,}\\
u_i&\mbox{with probability $\eta$,}
\end{array}\right.
\end{equation*}
where $0\leq \eta<1$ denotes the (relative) noise level and $u_i$ is a number chosen randomly in the range of values of $\bfA \bfx_{\rm true}$. Although the images were generated using $\bfA(2.5,2.5)$, we consider reconstruction methods that use a different model matrix $\bfA(2.5,3.2)$, i.e., we introduce errors in the Gaussian blur parameters to model the realistic situation where the forward operator contains uncertainty and does not match data from the actual model.
We are interested in computing a value of $p$ such that the fit-to-data ${\rm L}^p$ norm in \eqref{eq:oid_inner} is suitable for handling the model mismatch.

Using the \OID\ approaches described in Section \ref{sub:learningexponents} with MM-GKS solvers for the inner problem, we compute the following optimal parameters.
\begin{itemize}
\item First, for fixed $p$ and $q$ we learn $\lambda$ only. For example, we denote `OID$_{\lambda,2,2}$' as OID with $\theta=\lambda$ and $p=q=2$, and we denote `OID$_{\lambda,1,2}$' as OID with $\theta=\lambda$, $p=1$, and $q=2$. The value of the regularization parameters so obtained are $0.4040$ and $0.4005$, respectively.
\item Then we learn the regularization parameter, fit-to-data norm, and regularization norm triplet.  We refer to this approach as `OID$_{\lambda,p,q}$' and we obtained the values of $\widehat\bftheta=[\ \widehat\lambda,\ \widehat p,\ \widehat q \ ]\t= [0.3427, 0.8129, 0.4795]\t$.
\end{itemize}
All OID approaches use surrogate optimization with a maximum of $200$ iterations and with lower and upper bounds of $ 10^{-8} \leq \lambda \leq 10$ for OID$_{\lambda,2,2}$ and OID$_{\lambda,1,2}$, and lower and upper bounds of $ 10^{-8} \leq \lambda \leq 10$, $ 0.1 \leq p, q \leq 2.5$ for OID$_{\lambda,p,q}$.  For the inner problem, we prescribed $50$ iterations of MM-GKS and $100$ iterations of CGLS (for $p=q=2$) at both training and validation stages. For all of the results in this section, we use the sample mean to center the data, prior to learning.

\begin{figure}
    \centering
    \includegraphics[width = 0.65\textwidth]{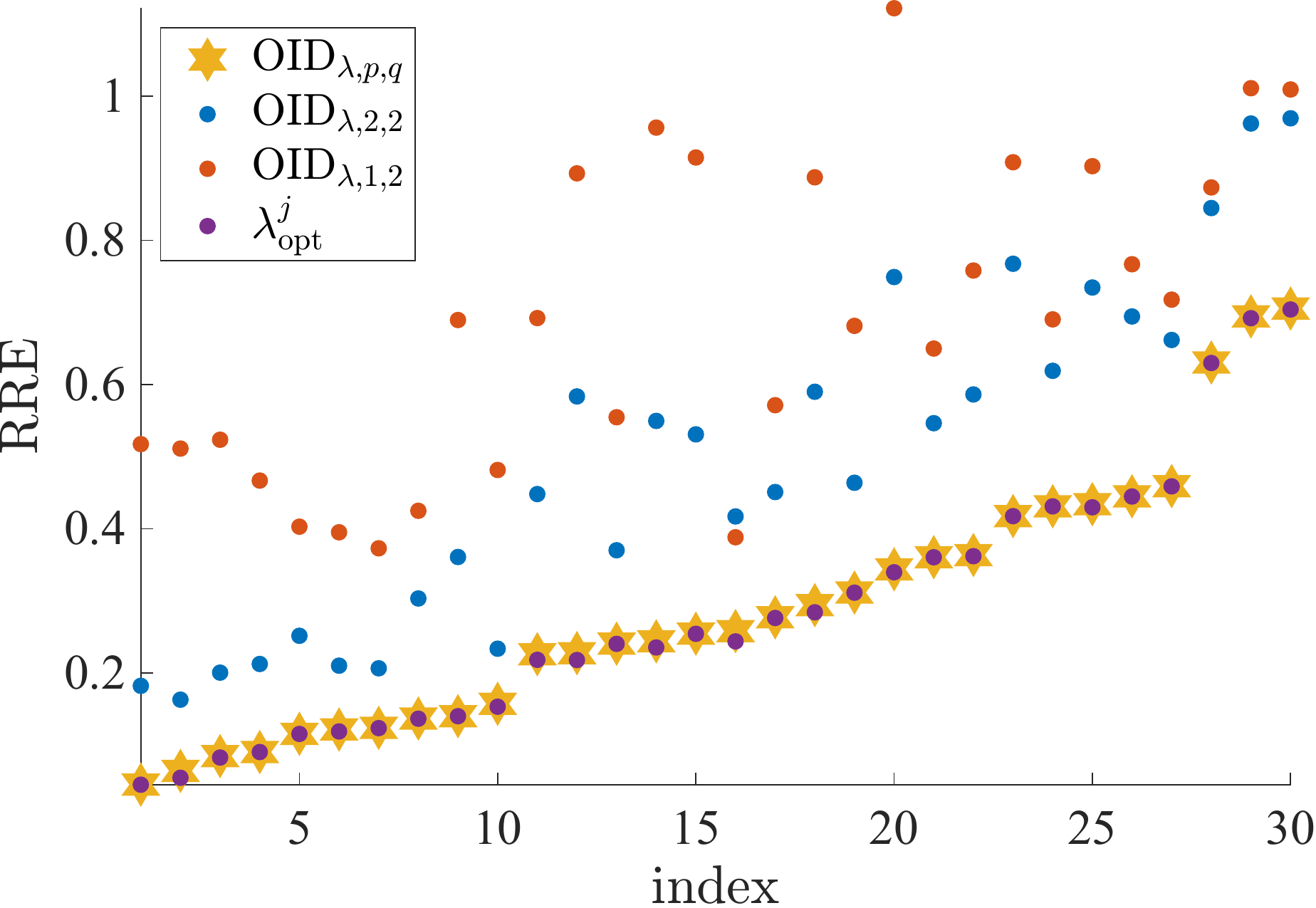}
    \caption{For the image deblurring problem with model error and impulse noise, we provide scatter plots of RRE norms for OID$_{\lambda,p,q}$, OID$_{\lambda,2,2}$ and OID$_{\lambda,1,2}$.
    Each column of dots corresponds to one sample from the validation set, where the indices have been sorted based on the RRE norms for OID$_{\lambda,p,q}$.
    As a further comparison, $\lambda_{\rm opt}^j$ corresponds to RRE norms for~\eqref{traincompare}, where the optimal regularization parameter is selected for each image using the learned $\widehat p$ and $\widehat q$.}
    \label{fig:RREsat}
\end{figure}

\begin{figure}
    \centering
        \begin{tabular}{ccc}
        $\bfx_{\rm true}$ & $\bfb$ (RRE = 0.4402) & OID$_{\lambda,p,q}$ (RRE = 0.1238) \\
        \includegraphics[width=0.32\textwidth]{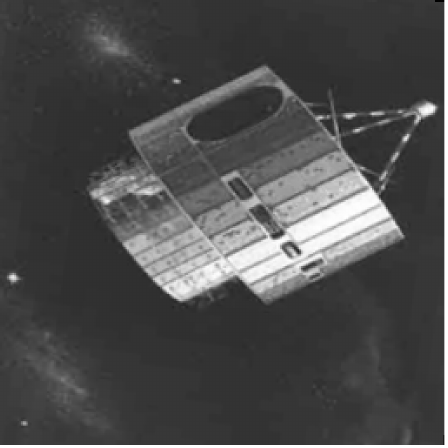} &
        \includegraphics[width=0.32\textwidth]{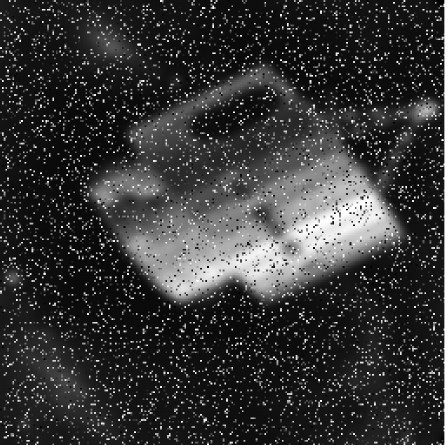}  &
        \includegraphics[width=0.32\textwidth]{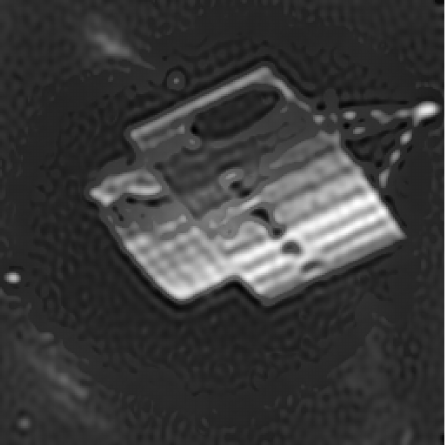} \\
        OID$_{\rm \lambda,2,2}$ (RRE = 0.2064) & OID$_{\lambda,1,2}$ (RRE = 0.3728) & ${\lambda_{\rm opt}^j}$  (RRE = 0.1235)\\
       \includegraphics[width=0.32\textwidth]{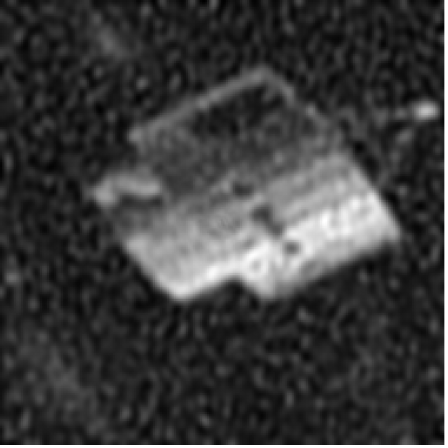} &
       \includegraphics[width=0.32\textwidth]{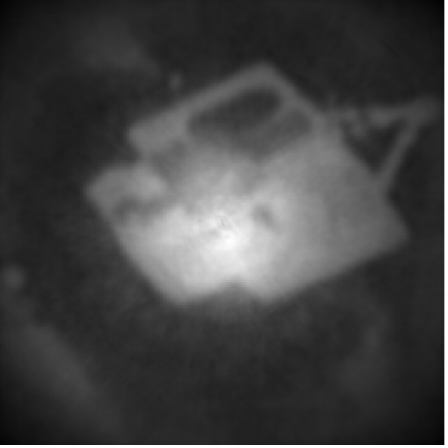} &
       \includegraphics[width=0.32\textwidth]{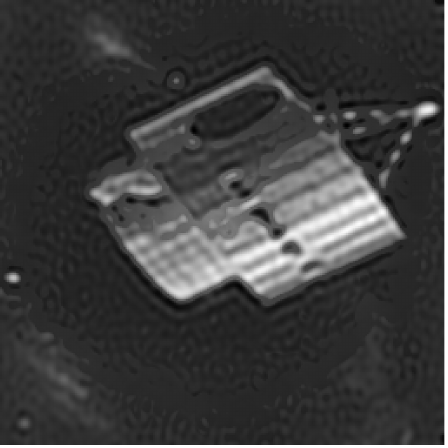} \\
    \end{tabular}
    \caption{For one sample of the validation data set, we provide in the top row the true image, the observed image, and the OID$_{\lambda,p,q}$ reconstruction. In the bottom row are reconstructions for OID$_{\lambda,2,2}$ and OID$_{\lambda,1,2}$. The bottom right reconstruction was computed using the optimal regularization parameter for this image and is provided for comparison purposes only.  RRE values are provided in the titles.}
    \label{fig:reconstrsat}
\end{figure}

Using OID computed parameters, we obtain reconstructions for each of the validation images. In Figure~\ref{fig:RREsat} we provide the RRE norms for OID$_{\lambda,p,q}$ (marked by yellow stars) for each validation image, where the index for the validation set has been sorted based on the RRE norms for OID$_{\lambda,p,q}$. RRE norms for OID$_{\lambda,2,2}$ (blue dots) and OID$_{\lambda,1,2}$ (red dots) are also provided for each validation image.
Since most of the blue and red dots lie above the yellow stars, we conclude that OID$_{\lambda,p,q}$ consistently performs better than OID$_{\lambda,2,2}$ and OID$_{\lambda,1,2}$, as expected. We also notice that RREs for OID$_{\lambda,2,2}$ are often smaller than RREs for OID$_{\lambda,1,2}$. Nevertheless, we observe that OID$_{\lambda,1,2}$ reconstructions eliminate the impulse noise well while struggling mitigate model errors in comparison to OID$_{\lambda,2,2}$ reconstructions, see Figure~\ref{fig:reconstrsat} (bottom left and middle).

To investigate the impact of the regularization parameter $\lambda$, we also provide results for the OID$_{\lambda,\widehat p,\widehat q}$ method, i.e., OID where the previously computed values for $p$ and $q$ are fixed.
Namely, the optimal regularization parameter $\lambda_{\rm opt}^j$ is computed for each validation image as
\begin{equation}
\begin{gathered}\label{traincompare}
\lambda_{\rm opt}^{j}=\argmin_{\lambda}\;\tfrac{1}{2} \norm[2]{\widehat\bfx^j(\lambda) -  \bfx_{\rm true}^j}^2\,,\\
\text{ where }\; \widehat{\bfx}^j(\lambda) = \argmin_{\bfx} \ \norm[\widehat p]{\bfA\bfx - \bfb^j}^{\widehat p} + \lambda \norm[\widehat q]{\bfx}^{\widehat q}\,,
\end{gathered}
\end{equation}
with the OID$_{\lambda,p,q}$ computed values $\widehat p=0.8129$, $\widehat q=0.4795$. RRE values for each validation image are provided as purple dots in Figure~\ref{fig:RREsat}. As expected, the images reconstructed using the optimal regularization parameter for each validation image have smaller RRE values than the images reconstructed using OID$_{\lambda,p,q}$ parameters.  However, we stress that this approach is not feasible in practice and that the OID results are not far off. Reconstructed images along with RRE values for one validation image are provided in Figure~\ref{fig:reconstrsat}. We observe that the OID$_{\lambda,p,q}$ reconstruction does not contain artifacts that are present in the OID$_{\lambda,2,2}$ and OID$_{\lambda,1,2}$ reconstructions (due to the learning of $p$ and $q$), and the reconstruction with the optimal regularization parameter is only slightly better and nearly indistinguishable from the OID$_{\lambda,p,q}$ reconstruction.

We remark that we observed similar performance for an example where the amount of blur is underestimated at the reconstruction stage; that is, training and validation data were generated using $\bfA(2.5,3.2),$ but $\bfA(2.5,2.5)$ was used for reconstructions. However, we do not provide results here.

\begin{figure}
    \centering\footnotesize
    \begin{tabular}{ccc}\small
    $\bfA(2.5,2.5)\bfx_{\rm true} - \bfb$ & $(\bfA(5.6,5.6) - \bfA(2.5,2.5))\bfx_{\rm true}$ & $\bfA(5.6,5.6)\bfx_{\rm true}- \bfb$\\
    \includegraphics[width = 0.32\textwidth]{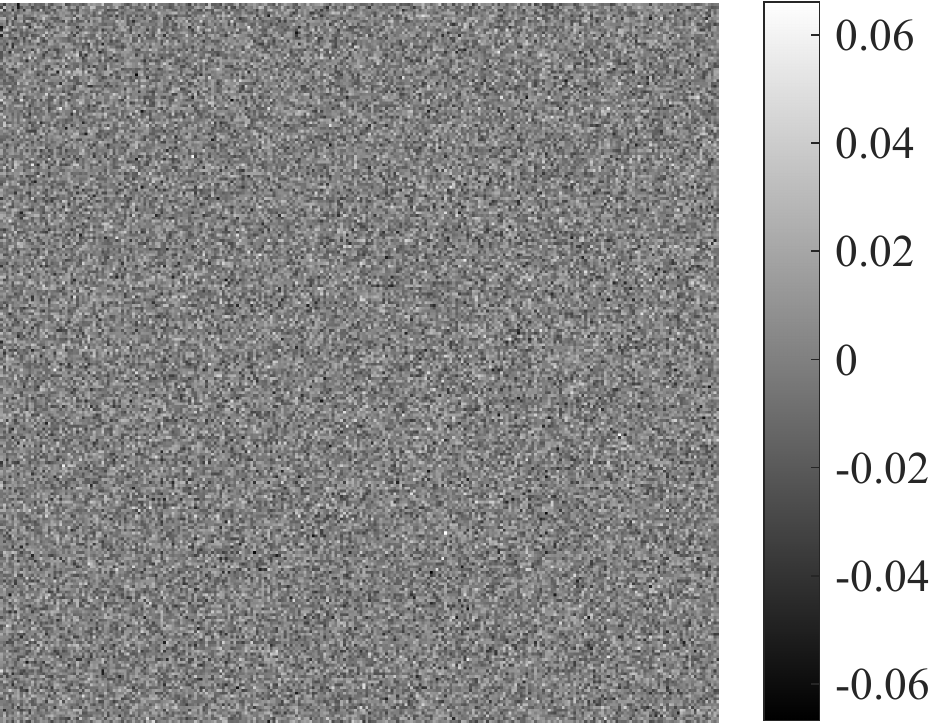} &
    \includegraphics[width = 0.32\textwidth]{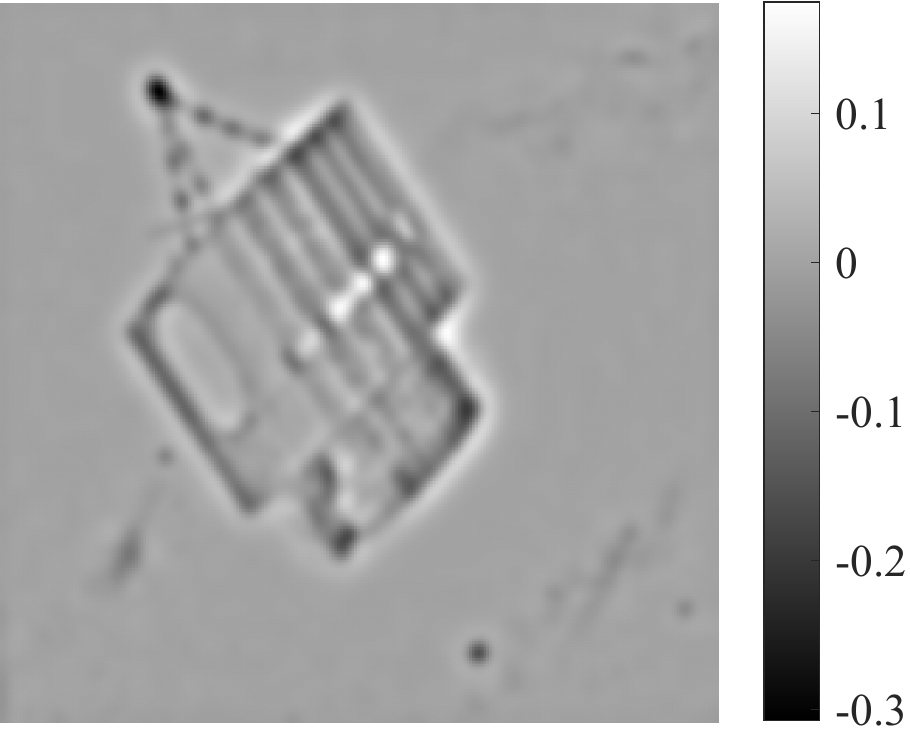} &
    \includegraphics[width = 0.32\textwidth]{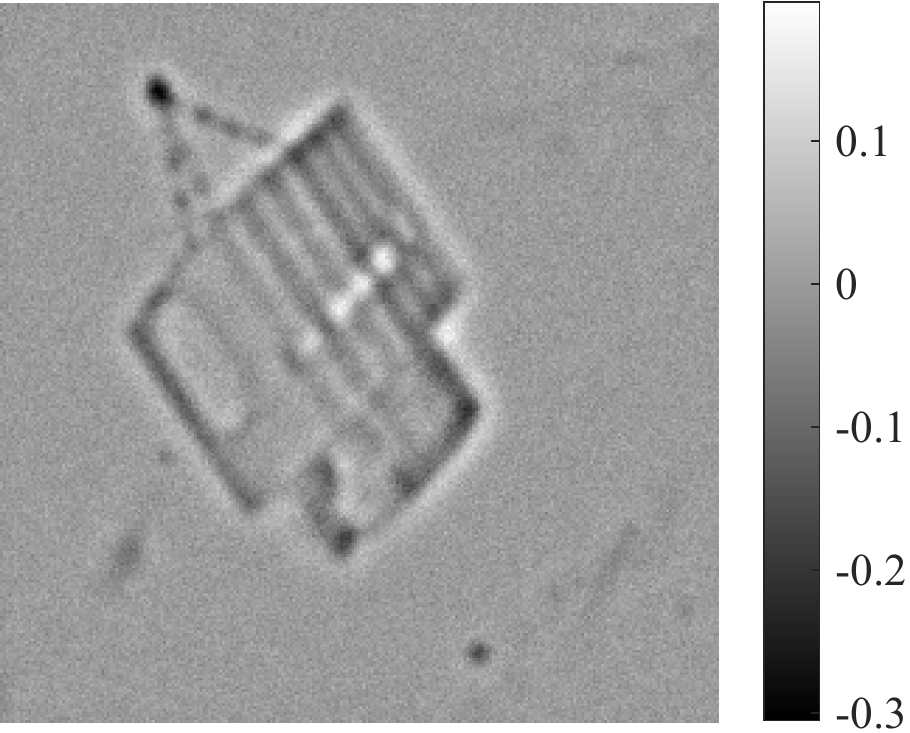}
\end{tabular}
    \includegraphics[width=.97\textwidth]{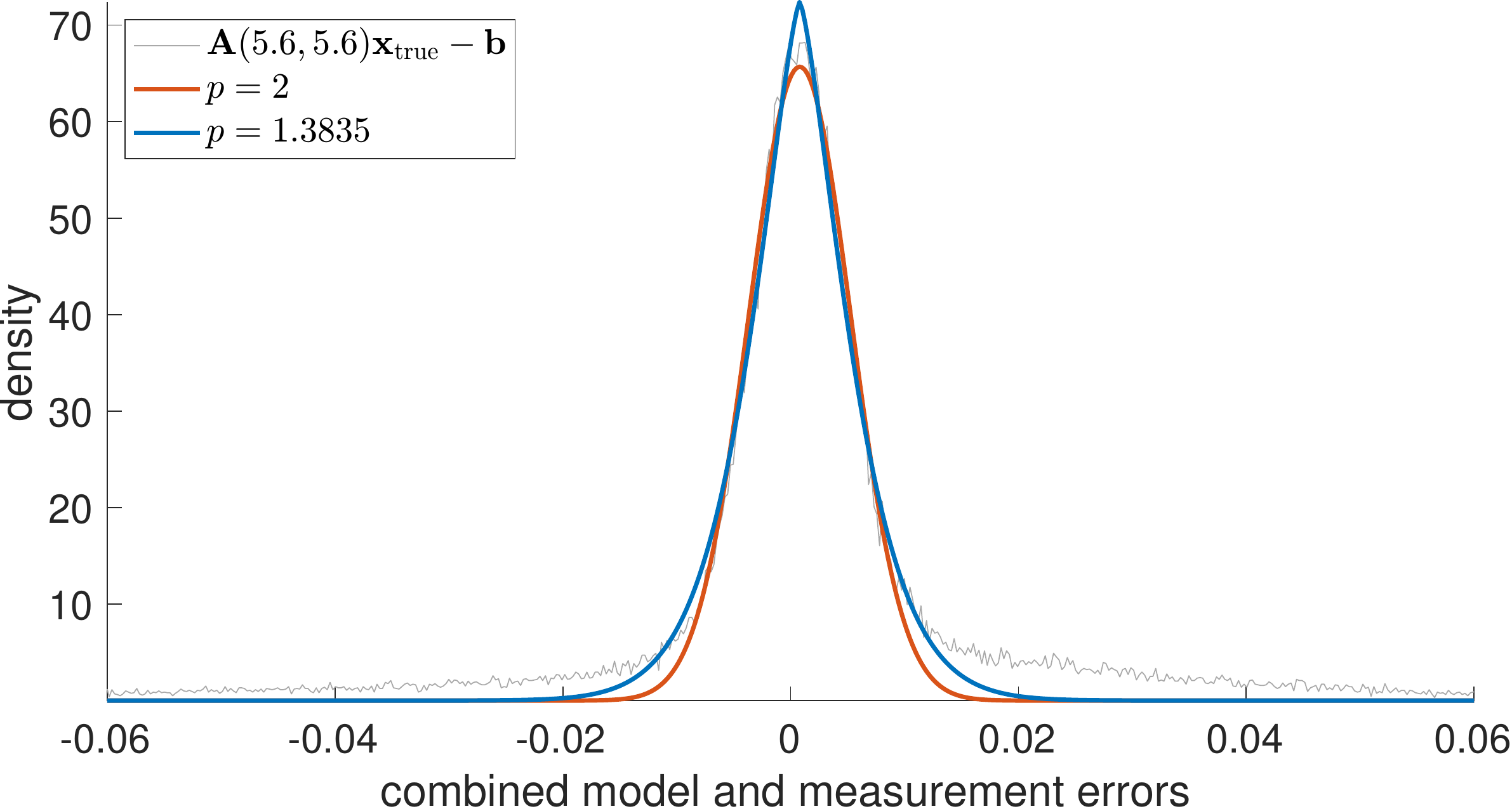}
    \caption{Investigating the impact of model error on the overall noise statistics. From left to right in the top row, we provide pixel-wise values of the additive Gaussian noise (of level 1\%), the model error, and the sum of these two errors. The density plot for the combined error is provided, along with the density function for $p=2$ (corresponding to Gaussian noise) and the best density fit to the true errors $\bfA(5.6,5.6)\bfx_{\rm true}- \bfb$, i.e., $p = 1.3835$.}
    \label{fig:satstats}
\end{figure}

Next, we investigate the impact of model error on the noise and hence the data-fit term, for which it is well-known that the choice of $p$ is directly related to the statistics of the observation error. For one satellite image $\bfx_{\rm true}$, we consider the observed image that is generated using $\bfA(2.5,2.5)$ and we consider observation errors coming from two sources: 1\% additive Gaussian noise and model error by using $\bfA(5.6,5.6)$ for reconstructions instead of $\bfA(2.5,2.5)$. Images of the additive Gaussian noise, the model error, i.e., $\bfA(5.6,5.6)\bfx_{\rm true} -\bfA(2.5,2.5)\bfx_{\rm true}$, and the sum of the two sources of errors are provided in the top row of Figure~\ref{fig:satstats} from left to right.  These represent pixel-wise quantities. In the lower frame, we provide a density plot of the combined error, along with the density functions corresponding to $p=2$ and $p=1.3835$ (the best $p$-norm density fit for this image).

We observe that the combined model and measurement error, which resembles a heavy-tailed distribution, results in a noise distribution that is no longer Gaussian (i.e., ignoring the model error and using $p=2$ may not be appropriate). Indeed, even with additive Gaussian noise, there exists a value for $p$ that better resembles the noise statistics when model error is present. Thus, without additional prior knowledge, changing the norm for the data-fit term (in effect learning the statistics of the combined additive and model error from training data) is a reasonable approach to handle model error.

\subsection{OID for learning kernel parameter(s) and regularization parameter}
\label{sub:numerical_kernel}
\begin{figure}[htb!]
		\centering
		\includegraphics[width = .45\textwidth]{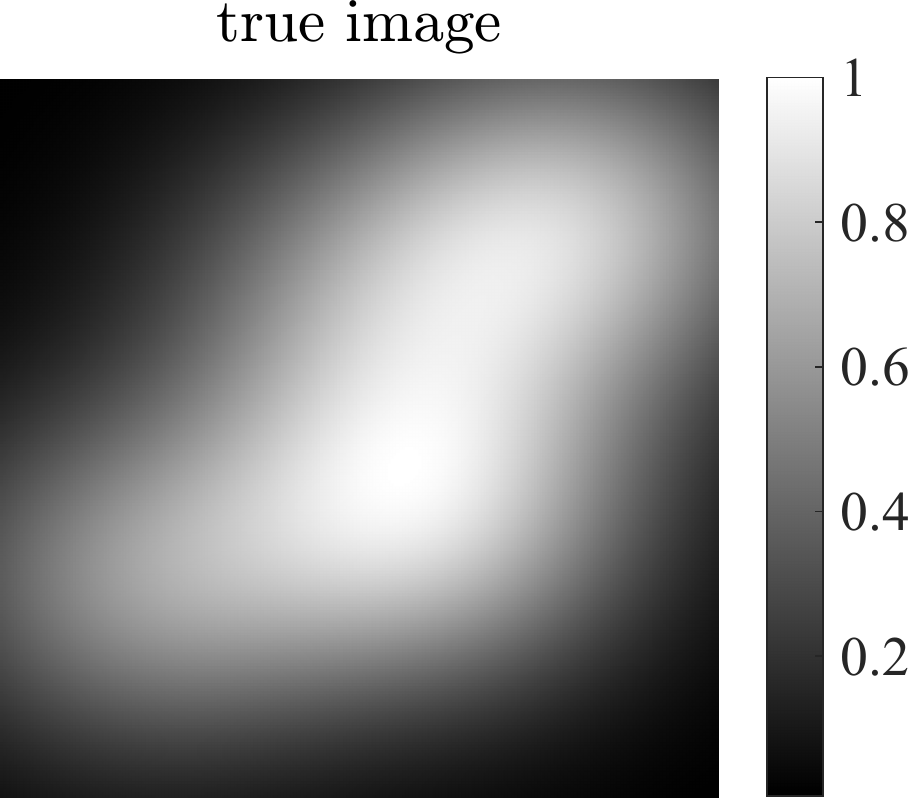}
		\includegraphics[width = .35\textwidth]{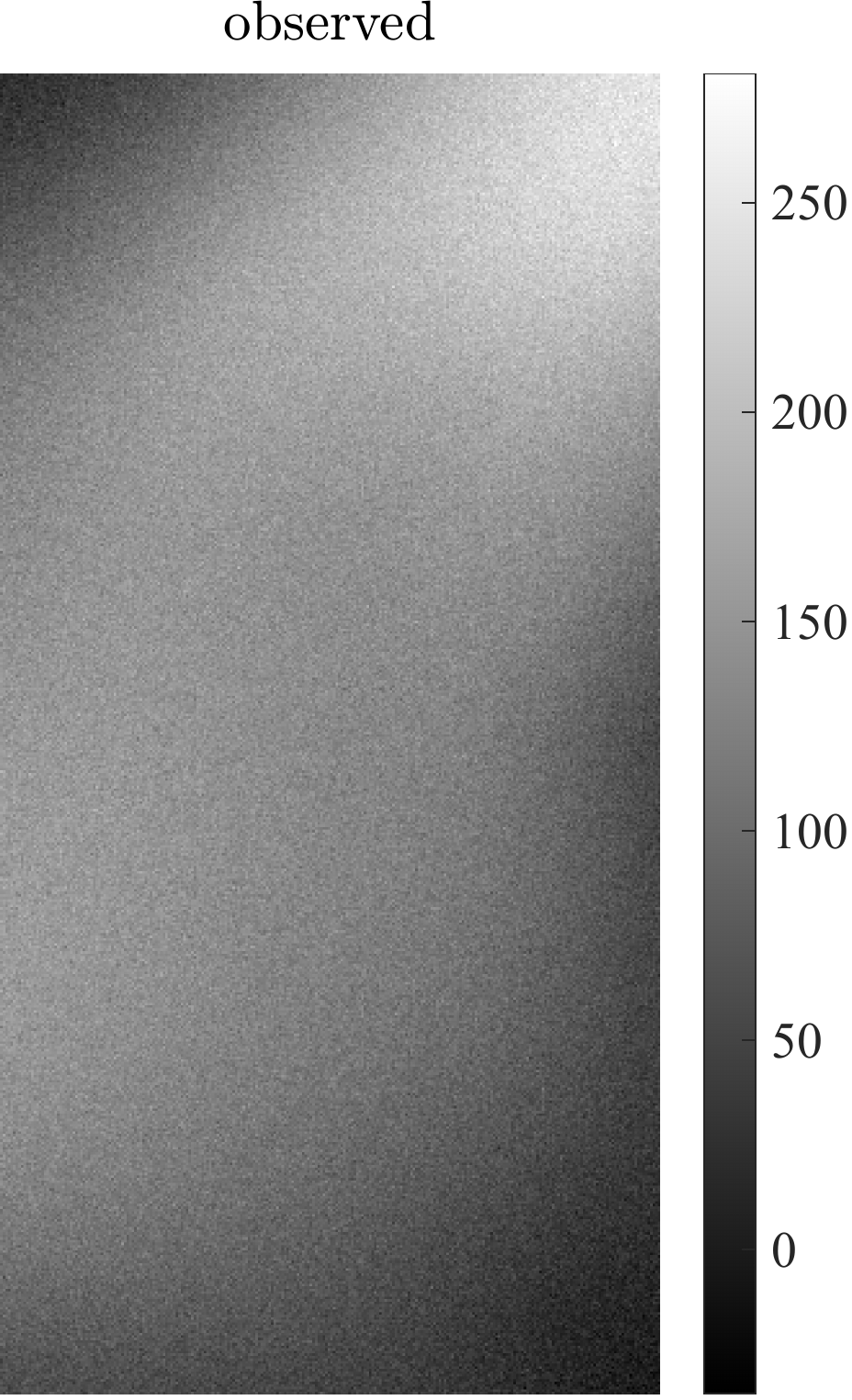}
\caption{Seismic image reconstruction example.  The true image (left) contains $256 \times 256$ pixels and represents a smooth medium.  The noisy sinogram image (right) represents projection data from a setup with $256$ sources and $512$ receivers.}
\label{fig:seismic_problem}
\end{figure}

We consider a seismic imaging problem (namely, PRseismic from  \cite{hansen2018air}) that can be modeled as \eqref{eq:linEq}, with $\bfx_{\rm true}$ containing a smooth image; see Figure~\ref{fig:seismic_problem}. $\bfA$ represents 2D seismic travel-time tomography, using $n_{\rm s}=256$ sources located on the right boundary and $n_{\rm r}=512$ receivers (seismographs) scattered along the left and top boundaries.  The rays are transmitted from each source to each receiver.  The noisy observations are provided in Figure~\ref{fig:seismic_problem};
here $n=65,\!536 = 256^2$.

\begin{figure}[h!]
		\centering
		\begin{tabular}{cccc}
		  \includegraphics[width = .22\textwidth]{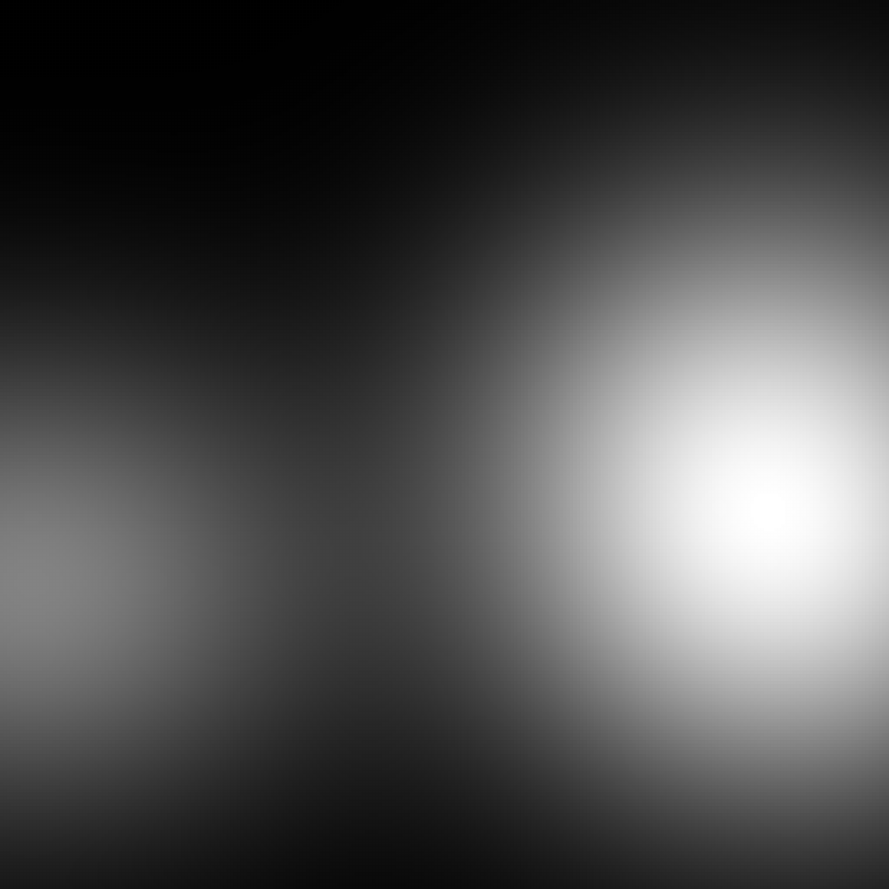} &
		  \includegraphics[width = .22\textwidth]{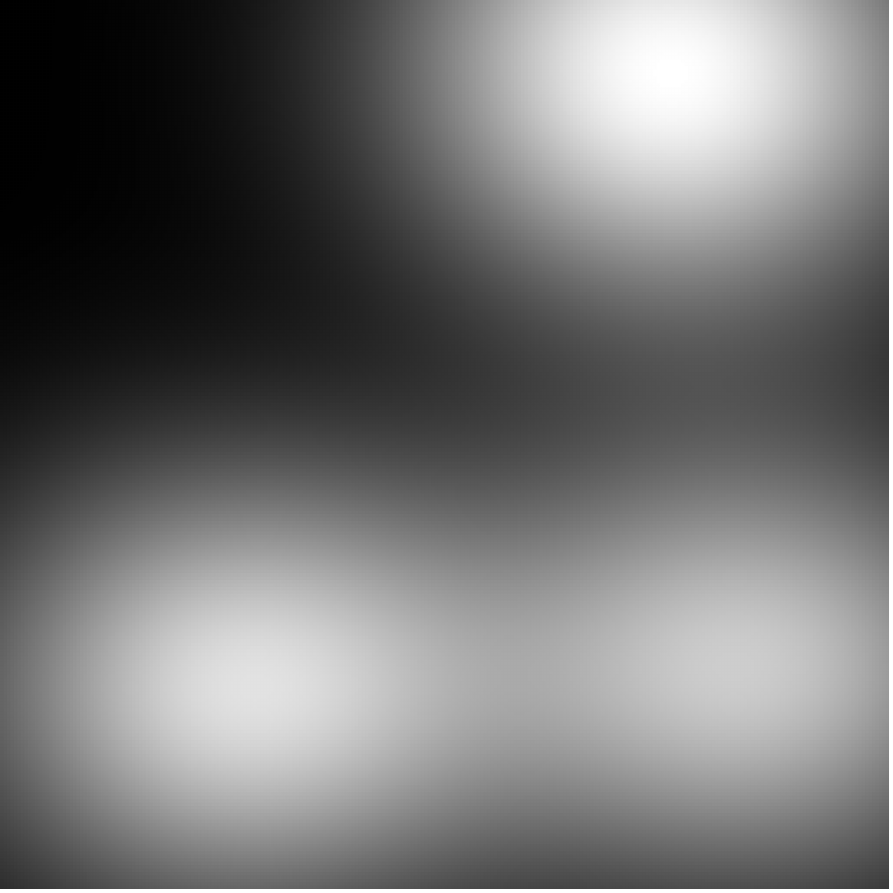} &
		  \includegraphics[width = .22\textwidth]{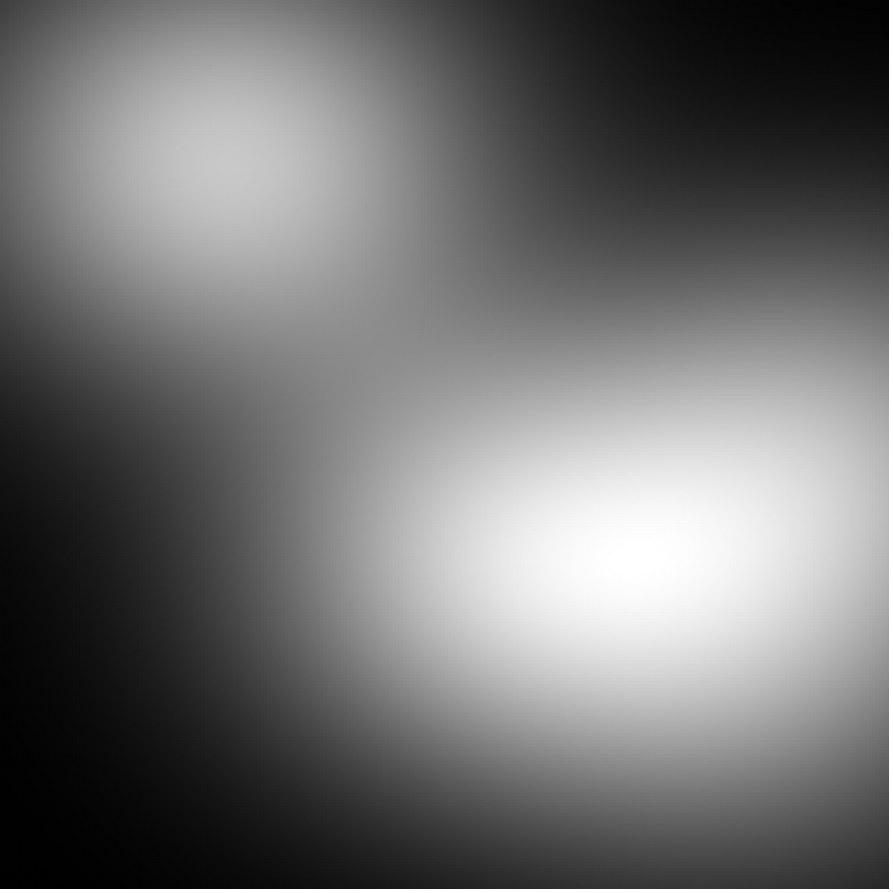} &
		  \includegraphics[width = .22\textwidth]{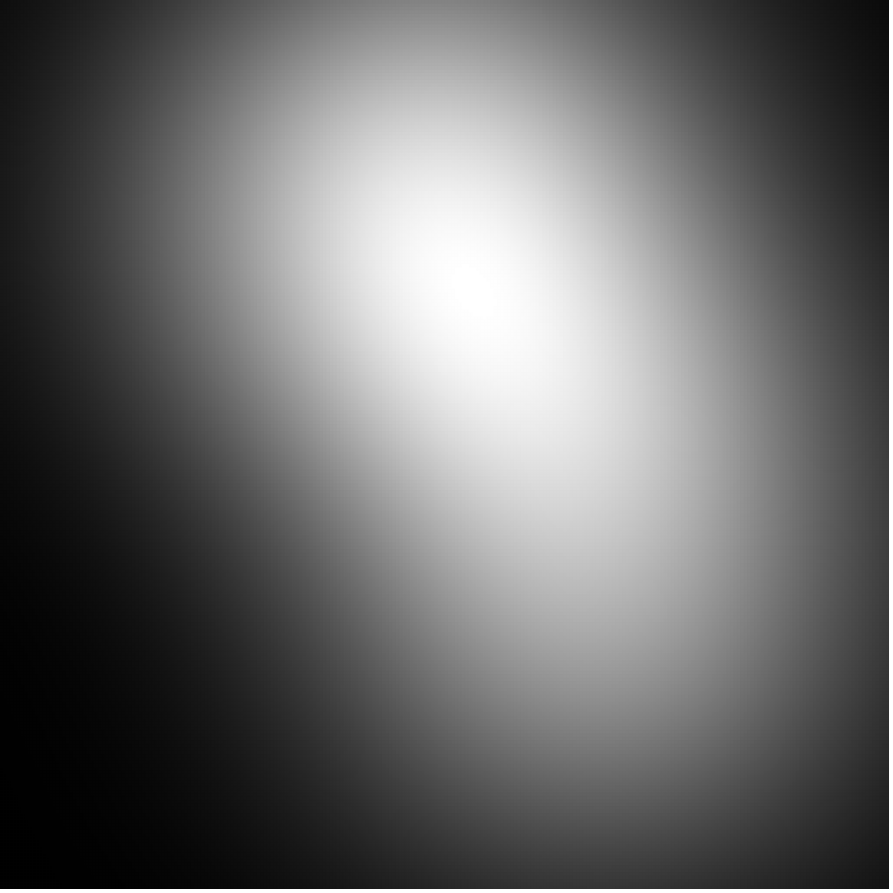} \\
		  \includegraphics[width = .22\textwidth]{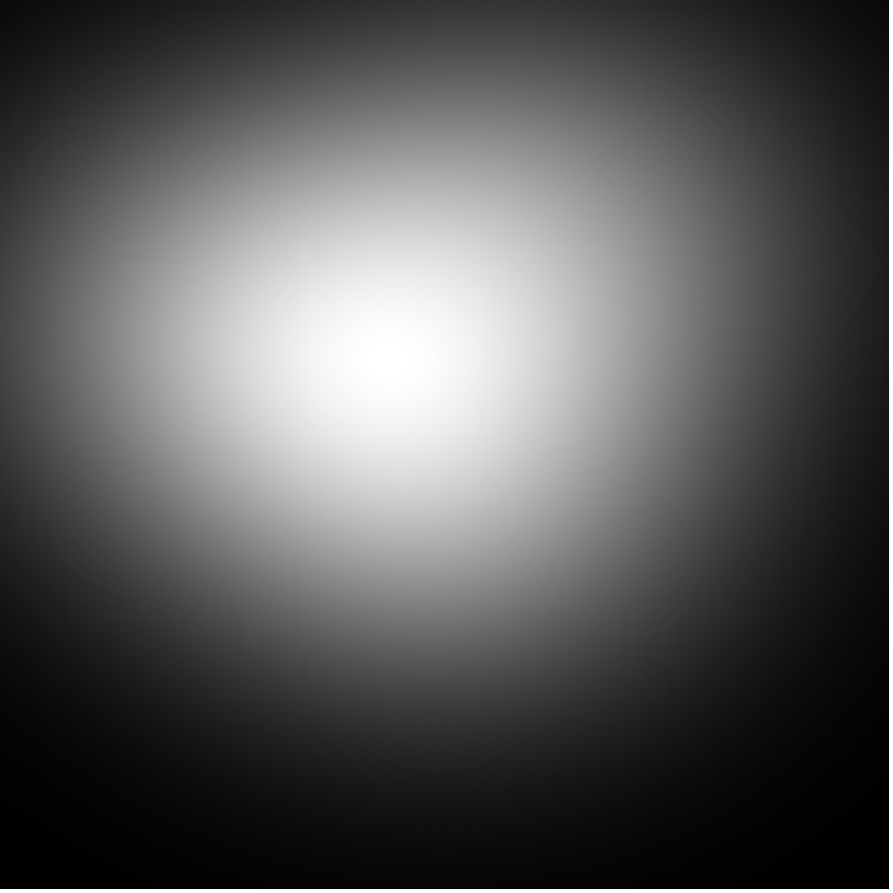} &
		  \includegraphics[width = .22\textwidth]{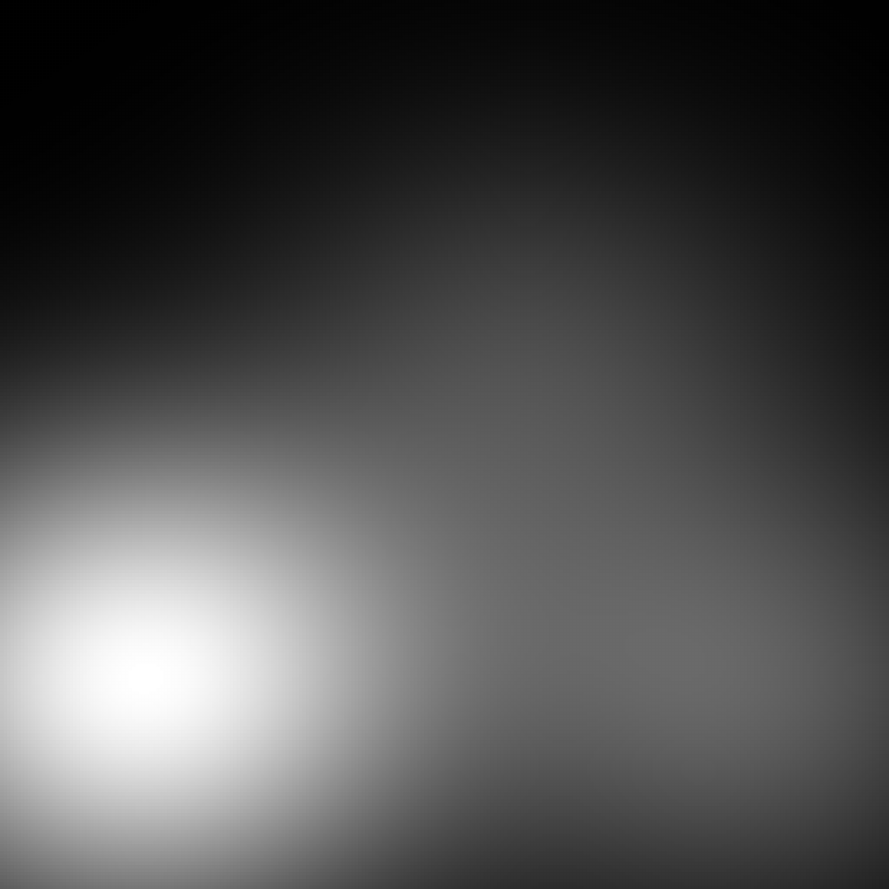} &
		  \includegraphics[width = .22\textwidth]{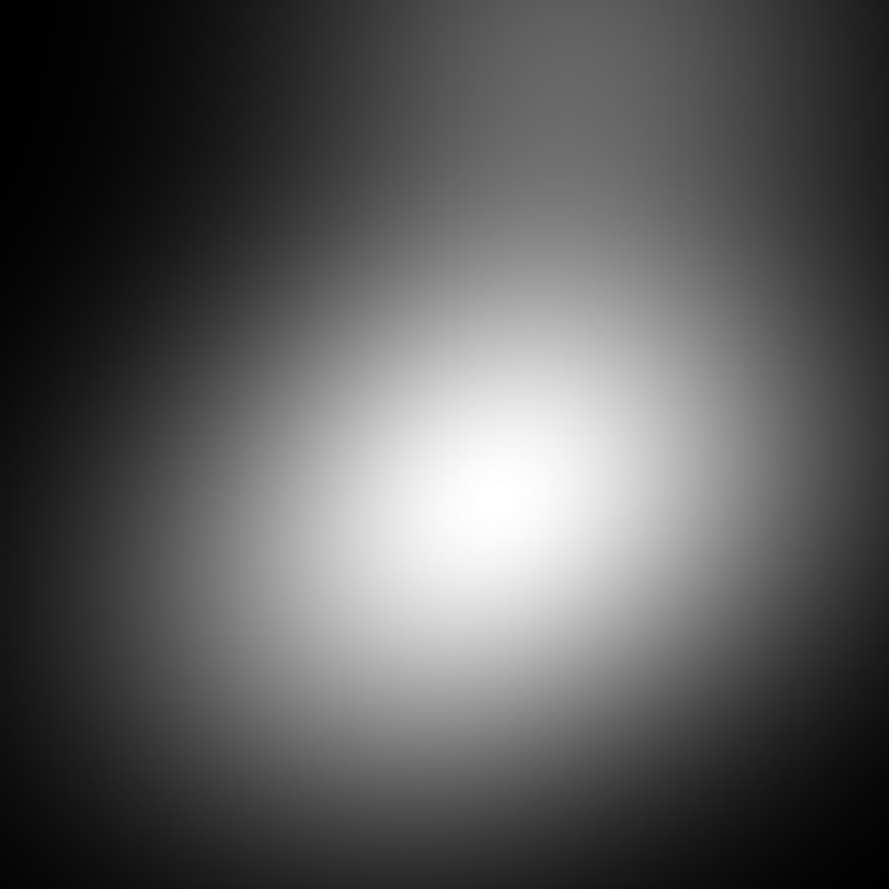} &
		  \includegraphics[width = .22\textwidth]{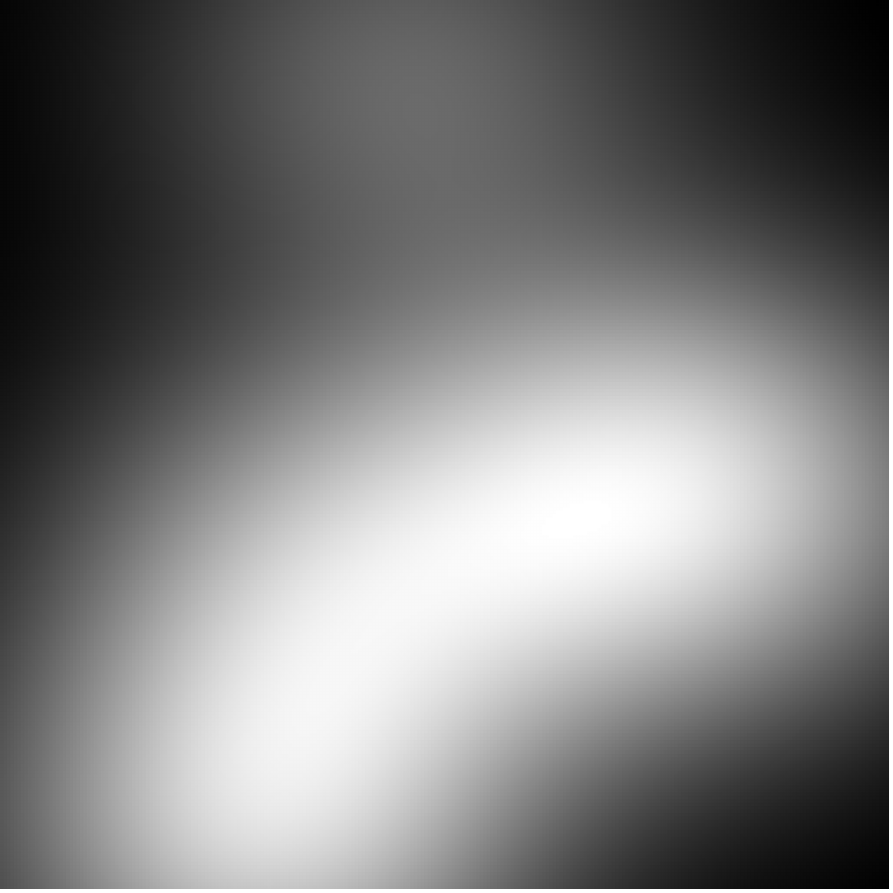}
		\end{tabular}
\caption{Seismic example - random samples from the training set}
\label{fig:seismic_samples}
\end{figure}

We generated a set of $30$ training images, some samples are provided in Figure~\ref{fig:seismic_samples}.  These were obtained by randomizing the parameters used to define the ``smooth'' image in IRTools \cite{gazzola2019ir}. Then for each training image, we simulated noisy observations using realizations of a Gaussian random vector with $\bfzero$ mean and noise level $\eta=\frac{\|\bfe\|_{2}}{\|\bfA\bfx_{\rm true}\|_{2}}$ uniformly chosen between $10^{-2}$ and $10^{-1}$.

Using the training data, we consider two kernel functions,
the squared exponential kernel function \eqref{eq:sqexp} and the Mat$\acute{\text{e}}$rn kernel function \eqref{eq:matern}. For each kernel function, we provide OID results for the following two scenarios:
\begin{itemize}
\item `OID' corresponds to solving \eqref{eq:oid} with $\bftheta = [\lambda, \bfbeta]$, where the genGK-based iterative projection method described in Section~\ref{sub:operators} is used to solve the inner problem \eqref{eq:oid_inner}.  Here we note that, since $\lambda$ is being learned in the OID problem, the only stopping criteria used for the inner problem are based on tolerances on the residual norm.
\item OID with $\bftheta = \bfbeta$, where genHyBR is used for selecting $\lambda$ according to WGCV (weighted generalized cross validation) and the full suite of stopping criteria are used within genHyBR. We refer to this approach as `OID-wgcv'.
\end{itemize}
Both OID approaches use surrogate optimization with a maximum of $20$ iterations and with lower and upper bounds of $10^{-6} \leq \lambda \leq 1$, $0.5 \leq \beta_1 \leq 15$, and $5\cdot 10^{-2} \leq \beta_2 \leq 0.7$ for the Mat$\acute{\text{e}}$rn kernel and $10^{-6} \leq \lambda \leq 1$ and $0.01 \leq \beta \leq 0.5$ for the squared exponential kernel. Computed \OID\ parameters are given in Table~\ref{tab:ComputedOID}.

\begin{figure}[h!]
		\centering
	  \includegraphics[width = .7\textwidth]{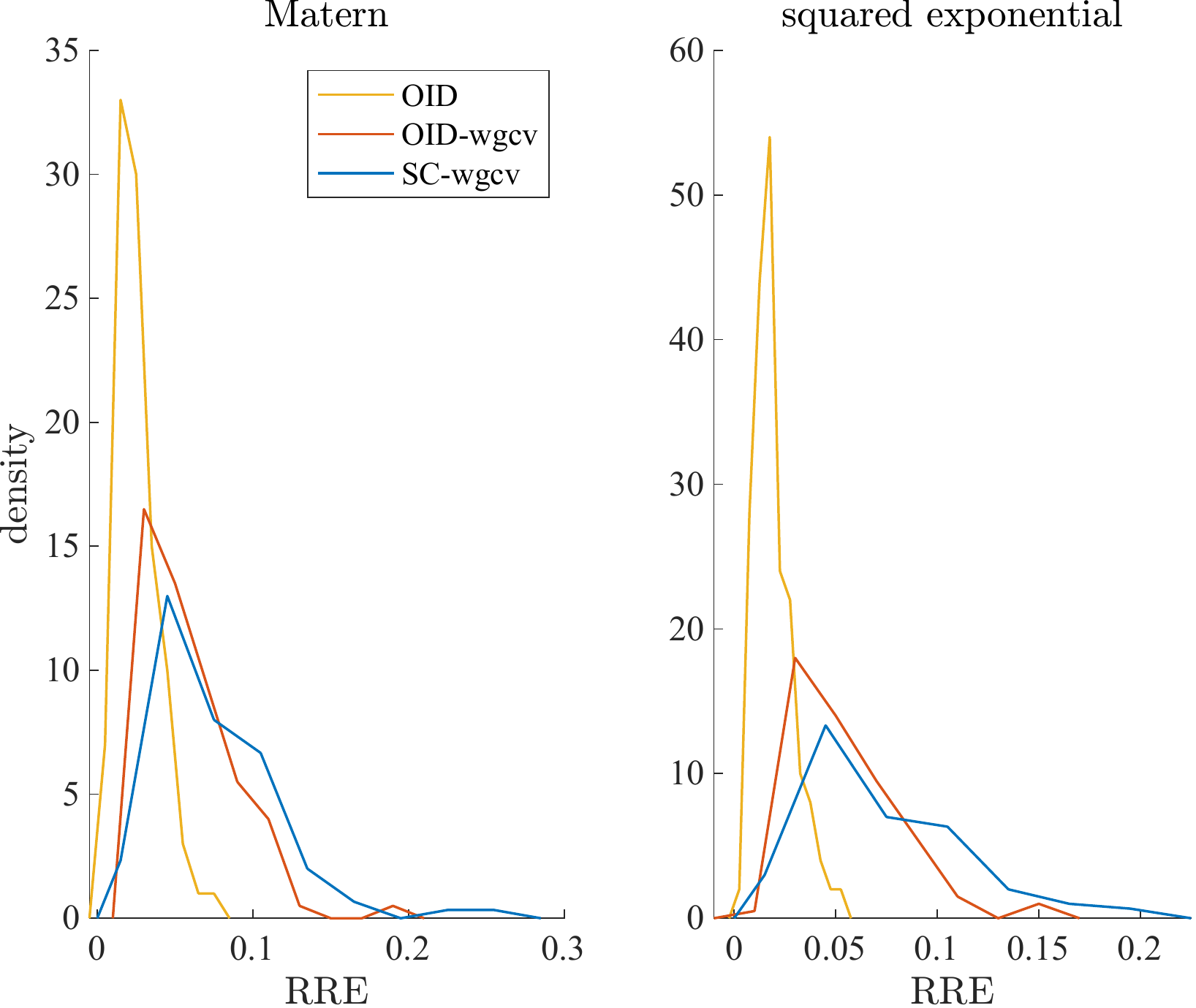}
\caption{For the validation set of the seismic example, we provide histograms of the RRE norms for OID, OID-wgcv, and SC-wgcv.}
\label{fig:seismic_hist}
\end{figure}

Then similar to how we generated the training images, we generated $100$ validation images and their corresponding observations.  Using the \OID\ parameters in Table~\ref{tab:ComputedOID}, we obtained reconstructions for the validation set. The overall mean squared reconstruction error for the validation images is provided in the last column of Table~\ref{tab:ComputedOID}.
Individually computed RRE norms for OID and OID-wgcv for each validation image are provided as yellow and red dots in Figure~\ref{fig:seismic_hist} respectively, where the indices have been sorted based on the RRE norms for OID.  Notice that most of the red dots lie above the yellow dots, and thus OID tends to perform better than OID-wgcv for this example.

\begin{figure}[htb!]
		\centering
	  \includegraphics[width = .9\textwidth]{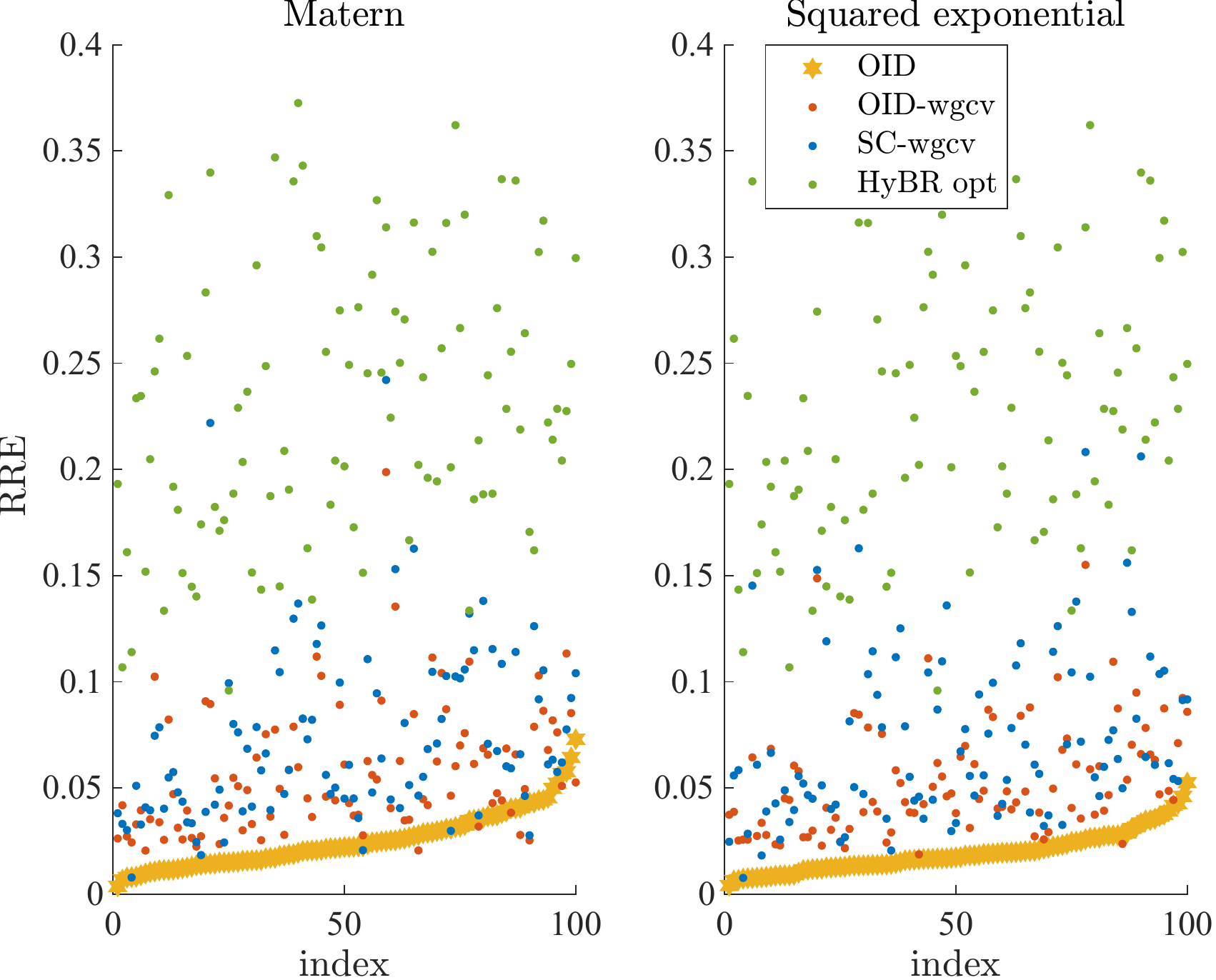}
\caption{For the seismic example, we provide scatter plots of RRE norms for OID, OID-wgcv, SC-wgcv, and HyBR opt. Each column of dots corresponds to one sample from the validation set, where the indices have been sorted based on the RRE norms for OID.}
\label{fig:seismic_scatter}
\end{figure}

\begin{figure}[h!]\footnotesize
\begin{tabular}{ccc}
 \hspace*{-7ex} $\lambda_{\rm opt}$ $\bfQ=\bfI$ & \hspace*{-7ex} OID Mat$\acute{\text{e}}$rn & \hspace*{-7ex} OID sq. exp.\\
    \includegraphics[width = .31\textwidth]{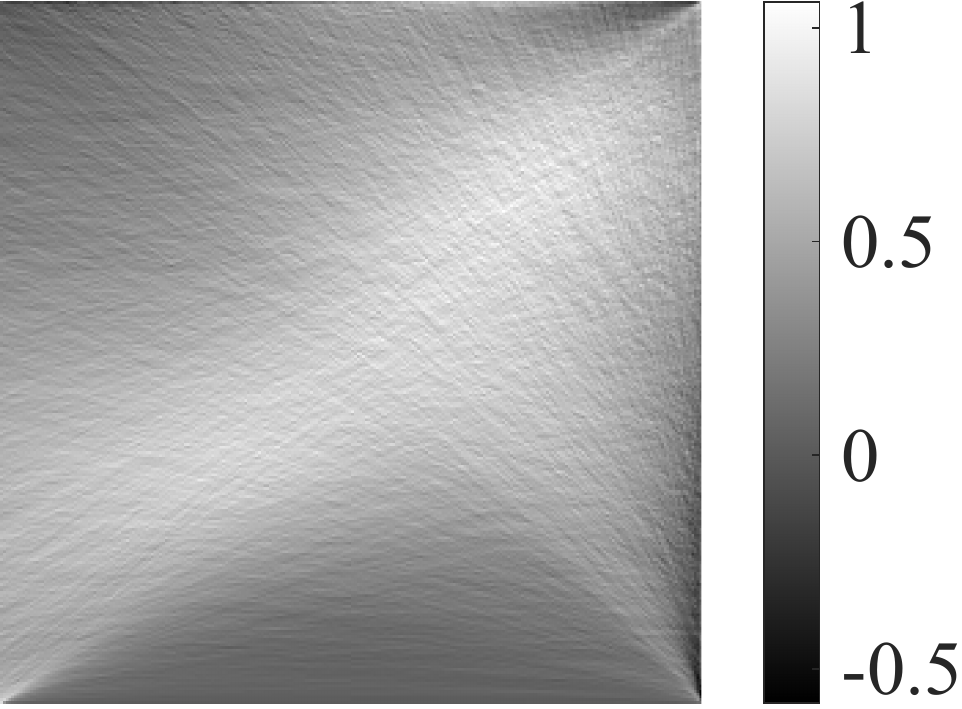} &
    \includegraphics[width = .31\textwidth]{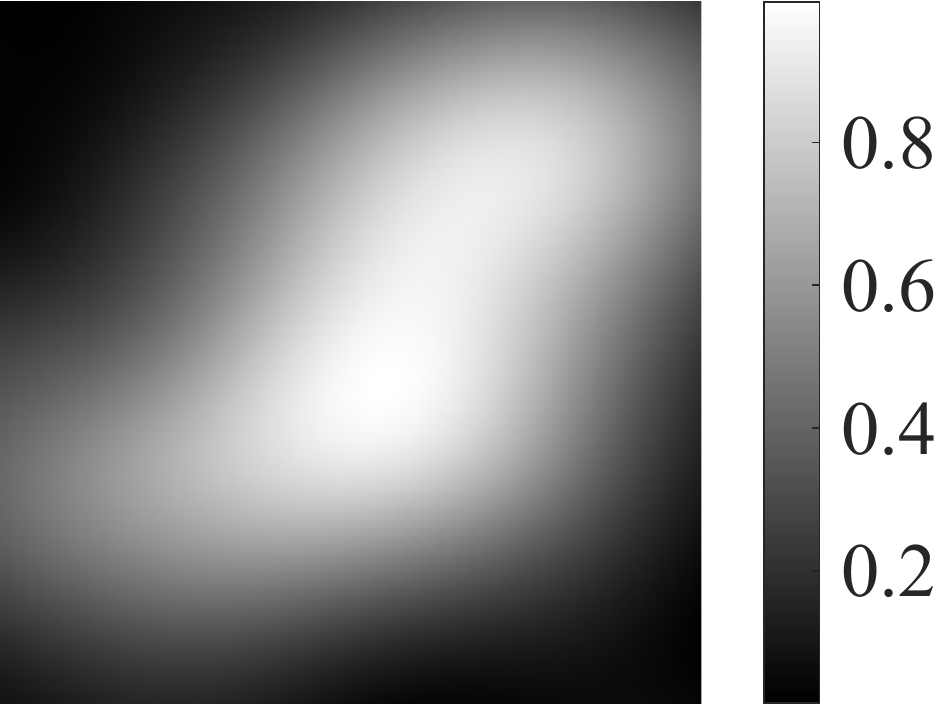} &
    \includegraphics[width = .31\textwidth]{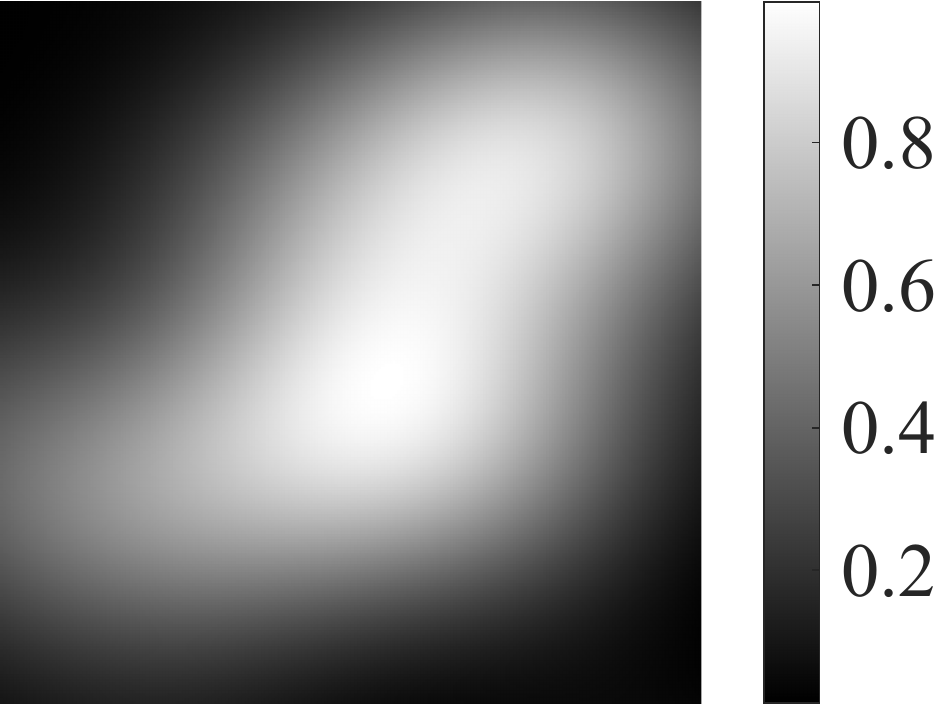}\\
    \hspace*{-7ex} (RRE = 0.2554) & \hspace*{-7ex}  (RRE = 0.0214) & \hspace*{-7ex}  (RRE = 0.0190)
\end{tabular}
\caption{For the validation image in Figure~\ref{fig:seismic_problem}, we provide reconstructions obtained with HyBR opt and OID reconstructions for the Mat$\acute{\text{e}}$rn and squared exponential kernels.}
\label{fig:seismic_val}
\end{figure}

For comparison, we use the approach described in \cite{cho2021hybrid} that estimates hyperparameters directly from the the sample covariance matrix constructed from the training data, followed by genHyBR with WGCV for selecting $\lambda$. Following \cite{cho2021hybrid}, let $\widehat\bfQ$ be the sample covariance matrix constructed from the training dataset, then Mat$\acute{\text{e}}$rn parameters were estimated by solving
\begin{equation}\label{Fapprox}
        (\hat \nu, \hat \ell) = \argmin_{\nu>0,\ell>0}\norm[\rm F]{\bfQ(\nu,\ell) - \widehat \bfQ}^2,
\end{equation}
where $\norm[\F]{\mdot}$ denotes the Frobenius norm.
Once the parameters are computed, they can be used to define $\bfQ = \bfQ(\hat\nu,\hat\ell)$, which can be used directly in generalized hybrid methods.
For computational feasibility, a stochastic approximation is used for the objective function in \eqref{Fapprox}, i.e.,
\begin{align}
    \norm[\F]{\bfQ(\nu,\ell) - \widehat{\bfQ}}^2
    & =  \bbE  \norm[2]{(\bfQ(\nu,\ell) - \widehat{\bfQ}) \bfxi}^2,
\end{align}
 where $\bfxi$ is a random variable such that $\bbE \bfxi = \bfzero$ and $\bbE (\bfxi \bfxi\t) = \bfI$.  Using a Hutchinson trace estimator, we let  $\bfxi^{(i)} \in \bbR^n$ for $i=1,2,\ldots,M$ be realizations of a Rademacher distribution (i.e., $\bfxi$ consists of $\pm1$ with equal probability), and we consider the approximate optimization problem,
 \begin{equation}
    \label{ref:RademacherEst}
    (\check{\nu},\check{\ell}) =\argmin_{\nu>0,\ell>0}
    \frac{1}{M}\sum_{i=1}^{M}  \|(\bfQ(\nu,\ell) - \widehat\bfQ)\bfxi^{(i)} \|^2_2\,.
\end{equation} Similar to the approach described in \cite{cho2021hybrid}, we used an interior-point method (\texttt{fmincon.m} in MATLAB) to minimize \eqref{ref:RademacherEst} with $M=100$, and refer to this approach as sample covariance (SC).  We extend this approach to be used for estimating the squared exponential kernel parameter $\beta$ and provide the computed hyperparameters in the row labeled `SC' in Table~\ref{tab:ComputedOID}.  We remark that this approach uses only the training data and not the model, noise or observations for learning the kernel parameter.  On the contrary, OID incorporates this information. Also, for comparison, we provide results for standard Tikhonov regularization ($\bfQ=\bfI$) with the optimal $\lambda$ selected for each sample.  The goal of this comparison is to show that including the prior is critical for this example.  Scatter plots of RRE values for both approaches are provided in Figure~\ref{fig:seismic_scatter}.  Density graphs of the reconstruction errors are provided in Figure~\ref{fig:seismic_hist}, and one reconstructed image is provided in Figure~\ref{fig:seismic_val}.
\begin{table}[bthp]
    \centering

    \begin{tabular}{|c|c|c|c|} \hline
    Mat$\acute{\text{e}}$rn & $\lambda$ & $\bfbeta$ & $\calP$, validation\\ \hline
OID & 18.8313 & 5.0312, 0.3344 & \textbf{6.0833} \\
OID & wgcv & 12.2812, 0.3344 & 29.9412 \\
SC & wgcv & 123.1735, 0.2011 & 49.3447 \\
\hline
sq. exp. & $\lambda$ & $\beta$ & $\calP$, validation\\ \hline
OID & 50.0500 & 0.2550 & \textbf{3.4468} \\
OID & wgcv & 0.3163 & 24.7547 \\
SC & wgcv & 0.2010 & 47.6977 \\ \hline
\end{tabular}

    \caption{Computed values of the hyperparameters for OID, along with the mean reconstruction errors for the validation set. OID with $\lambda$ computed using WGCV corresponds to using OID-wgcv for estimating $\bfbeta$ only. `SC' corresponds to estimating $\bfbeta$ directly from the sample covariance matrix as described in \cite{cho2021hybrid} and then using genHyBR with WGCV. }
    \label{tab:ComputedOID}
\end{table}

\begin{figure}[h!]
		\centering
	  \includegraphics[width = .8\textwidth]{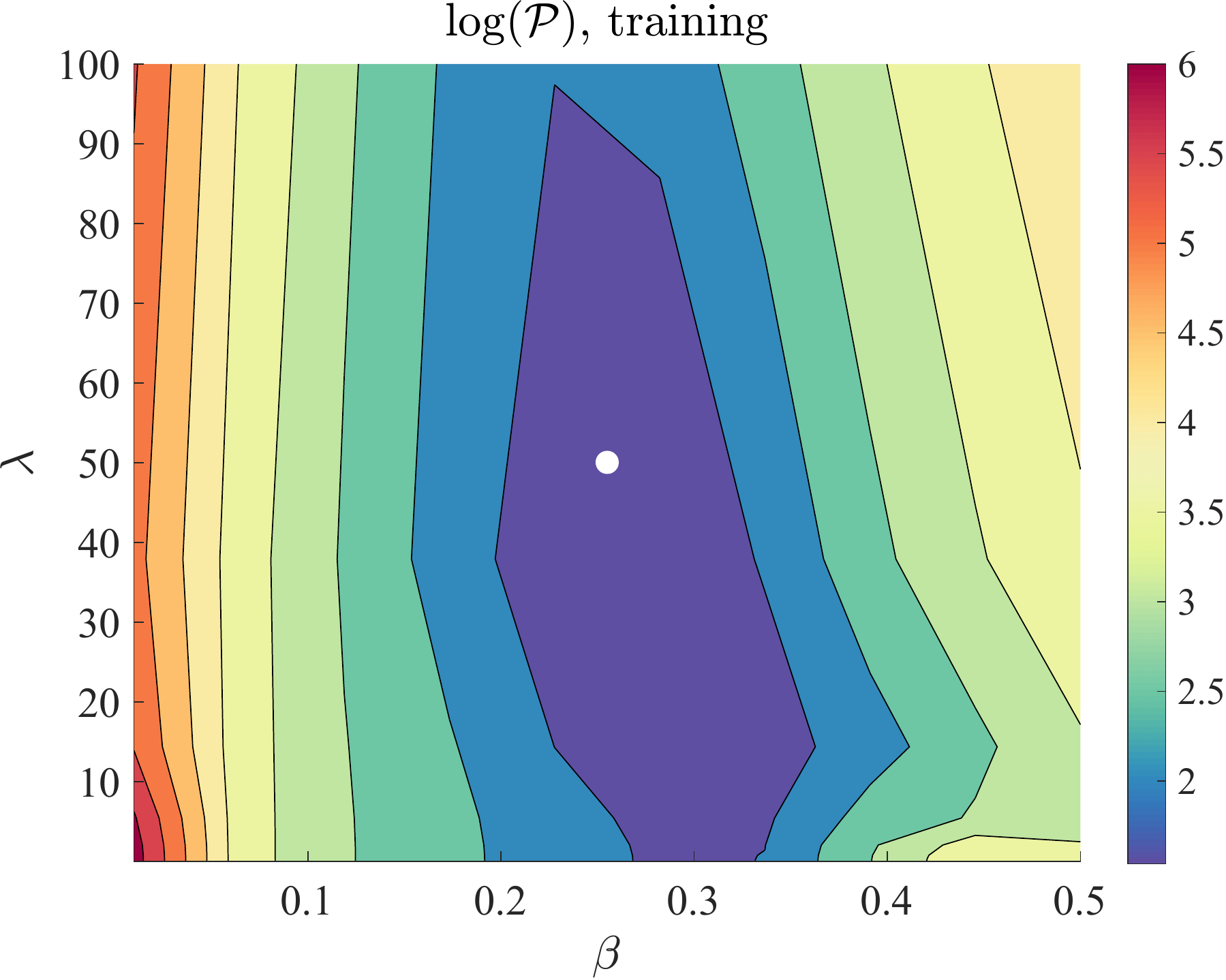}
\caption{Design objective for OID with the squared exponential kernel for the seismic example. The filled contour corresponds to OID, and the white point denotes the OID computed values.}
\label{fig:seismic_surface}
\end{figure}

Finally, we investigate the properties of the design objective function \eqref{eq:oid_outer} for OID with the squared exponential kernel.  In Figure~\ref{fig:seismic_surface}, we provide a contour plot of the design objective function, where the white point corresponds to the OID computed values.  Notice that there is wide region of values for $\lambda$ and $\beta$ that result in small and similar design objective values, with a smaller range of good choices for $\beta$.

\section{Conclusions and extensions}
\label{sec:conclusions}
In this work, we have presented a unified framework for optimal inversion design for large-scale inverse problems.  We have described learning approaches for computing hyperparameters from training data that exploit Krylov subspace methods for efficiently solving regularized inner problems within bi-level schema. In particular, we considered OID for learning the norm exponent in the data-fit and regularization term, as well as for learning the regularization parameter.  Furthermore, we considered OID for learning parameters for kernel functions used to define prior covariance matrices. Numerical experiments showed that OID methods can compute hyperparameters that deliver quality reconstructions, even in especially relevant scenarios where there is a mixture of noise and model error in the data (e.g., due to the presence of inaccuracies in the forward operator).

We remark that there are other cases where data-driven, optimal inverse frameworks can be used. The focus of our applications is image processing, and more particularly, image deblurring and computerized tomography; nevertheless, the learning approaches that we propose here can be applied to broad applications outside the field of image processing. Moreover, general regularization matrices can be used when considering OID with $\bftheta=[p,q,\lambda]\t$, including:  discretizations of the derivative operators when solutions with edge preserving properties are desired, wavelet and framelet transformations like in \cite{buccini2020modulus, BRP18, buccinibregman,COS09b,HDC13} when the solution is sparse in a transformed domain, or fractional Laplacian regularizers where smoothness is determined by a fractional exponent \cite{antil2020bilevel}. Other extensions include more general model errors, e.g., where the true forward model is a matrix perturbation of the forward model matrix, and stochastic approximation methods for problems where the training set is very large and hence empirical Bayes risk minimization problems become computationally intractable. These are topics of future investigations.

\section*{Acknowledgments}
This work was initiated as a part of the Statistical and Applied Mathematical Sciences Institute (SAMSI) Program on Numerical Analysis in Data Science in 2020. Any opinions, findings, and conclusions or recommendations expressed in this material are those of the authors and do not necessarily reflect the views of the National Science Foundation (NSF). MP gratefully acknowledges the  support from ASU Research Computing facilities for the computing resources used for testing purposes.
\printbibliography

\end{document}